\documentclass[12pt]{article}
\usepackage{amsmath,latexsym,amscd,amsfonts} 
\newcommand{\be}{\begin{equation}}
\newcommand{\ee}{\end{equation}}
\newcommand{\bes}{\begin{equation*}}
\newcommand{\ees}{\end{equation*}}

\numberwithin{equation}{section}
\newfont{\dr}{msym10 scaled \magstep1}
 \newfont{\gl}{eufm10 scaled \magstep1}  
 
 \newcommand{\bC}{{\bf C}}

 \newcommand{\bI}{{\bf I}}

 \newcommand{\bR}{{\bf R}}
 \newcommand{\bT}{{\bf T}}

 \newcommand{\cE}{{\cal E}}
 \newcommand{\cF}{{\cal F}}

 \newcommand{\cM}{{\cal M}}

 \newcommand{\lra}{\longrightarrow}
 
 \newcommand{\kahler}{K\"{a}hler}
 \newcommand{\noi}{\noindent}
 \newcommand{\dds}{\frac{\partial}{\partial s}}
 \newcommand{\ddt}{\frac{\partial}{\partial t}}

 \newcommand{\dbar}{\overline{\partial}}
 \newtheorem{definition}{Definition}[section]
 \newtheorem{lemma}[definition]{Lemma}
 \newtheorem{prop}[definition]{Proposition}
 \newtheorem{thm}[definition]{Theorem}
 \newtheorem{cor}[definition]{Corollary}

 \newcommand{\pf}{{\em Proof}. }
 \newcommand{\remark}{\noindent{\bf Remark}. }
 \newcommand{\remarks}{\noindent {\bf Remarks}. }

 \newcommand{\qed}{\hfill$\Box$}
 
 \newcommand{\kler}{\kahler}

\newcommand{\Hextn}{{0\lra (\cE_1,\Theta_1)\lra
 (\cE,\Theta)\lra(\cE_2,\Theta_2)\lra 0}}
\newcommand{\HopExtn}{{0\lra (E_1,\nabla_1'')\lra
 (E,\nabla'')\lra(E_2,\nabla_2'')\lra 0}} 
\newcommand{\Hsextn}{{0\lra(\cE_1',\Theta'_1)\lra(\cE',\Theta')\lra
(\cE_2',\Theta'_2)\lra 0}}

\newcommand{\dbare}{{\overline{\partial}_E}}
\newcommand{\delbar}{\overline{\partial}}
\newcommand{\lb}{\left(}
\newcommand{\rb}{\right)}
\newcommand{\Met}{\cM et}
\newcommand{\aHHE}{$\alpha$HHE}
\newcommand{\Fhiggs}{F^{\nabla}_H}
\newcommand{\SE}{\mathcal{E}}         
\newcommand{\SF}{\mathcal{F}}        
\newcommand{\EE}{\mathfrak{E}}   
\newcommand{\FF}{\mathfrak{F}}   
\newcommand{\PP}{\mathbb{P}}
\newcommand{\CC}{\mathbb{C}}
\newcommand{\SO}{\mathcal{O}}
\newcommand{\MH}{M_{Higgs}(X,P)}
\newcommand{\MHP}{M^\alpha_{Higgs}(X,P,P'')}
\newcommand{\MPP}{M^\alpha_{pure}(Z,\wt P,\wt P'')}
\newcommand{\MP}{M_{pure}(Z,\wt P)}
\newcommand{\MTF}{M^\alpha_{tf}(X,P,P'')}
\newcommand{\rk}{\operatorname{rk}}
\newcommand{\gr}{\operatorname{gr}}
\newcommand{\Gr}{\operatorname{Gr}}
\newcommand{\supp}{\operatorname{Supp}}
\newcommand{\slv}{\operatorname{SL(V)}}
\newcommand{\glv}{\operatorname{GL(V)}}

\newcommand{\quot}{\operatorname{Quot}}
\newcommand{\wt}{\widetilde}

\begin{document}
\begin{titlepage}
\title{Extensions of Higgs Bundles}
\author{Steven B. Bradlow \and Tom\'as L. G\'omez}
\thanks{The first author thanks the Tata Institute for Fundamental 
Research for its hospitality. The second author was supported by
a postdoctoral fellowship of Ministerio de Educaci\'on y Cultura
(Spain).} 
\maketitle
{\abstract We prove a Hitchin-Kobayashi correspondence for 
extensions of Higgs bundles. The results generalize known results 
for extensions of holomorphic bundles.  Using Simpson's methods, 
we construct moduli spaces of stable objects. In an appendix we 
construct Bott-Chern forms for Higgs bundles} 
 
\end{titlepage}
\section{Introduction}\label{sectn:intro}
\bigskip
The underlying principle at work in this paper is that, when 
approached in the right way, all results about holomorphic bundles 
can be made applicable to Higgs bundles too. 

\par 
The type of results we have in mind fall under the general heading 
of the Hitchin-Kobayashi Correspondence, i.e. they concern notions 
of stability, construction of moduli spaces, and the relation of 
these to solutions of gauge theoretic equations. Originally 
established for holomorphic bundles, results of this sort have 
been extended to Higgs bundles and also to a host of so-called 
`augmented holomorphic bundles', i.e. holomorphic bundles with 
some kind of prescribed additional structure.  Indeed a Higgs 
bundle can be treated as an augmented holomorphic bundle in which 
the augmentation is the Higgs field. However this is not always 
the best point of view - and is not the one we have in mind. The 
better approach is the one developed by Simpson in 
\cite{S1},\cite{S2},\cite{S3}. 

\par In Simpson's approach, instead of treating the Higgs structure 
as an augmentation, it is encoded in a more fundamental way. In 
fact there are two versions of this approach, one differential 
geometric and one algebraic. In the first (described in Section 
\ref{sectn:Eqtns}), the extra structure of a Higgs bundle is encoded as a 
modification of the partial differential operator which defines 
the holomorphic structure on a complex bundle. In the second (cf. 
Section \ref{sect:GIT}) , locally free coherent analytic sheaves 
on a variety $X$ are replaced by sheaves of pure dimension on 
$T^*X$. Having made 
these adjustments, a proof designed for holomorphic bundles or 
coherent analytic sheaves re-emerges as a proof for Higgs bundles 
or Higgs sheaves! 

\par In this paper we apply these principles to extensions of 
holomorphic bundles. A Hitchin-Kobayashi correspondence for such 
extensions was investigated in \cite{BGPcag} and \cite{DUW}; 
natural gauge-theoretic condition for special metrics, and a 
notion of stability were formulated, and the correspondence 
between them established. In \cite{DUW}, GIT methods were used to 
construct the moduli spaces. The main results in this paper thus 
show how, after the appropriate modifications,  these ideas can be 
carried over to Higgs bundles.  We set up and prove the 
Hitchin-Kobayashi correspondence for extensions of Higgs bundles 
(Theorems \ref{thm:easyHK} and \ref{thm:HardHK}), and we give (in 
Section \ref{sect:GIT}) a GIT construction for the associated 
moduli spaces. 

\par We also use the gauge-theoretic equations to deduce Bogomolov-type 
inequalities on the chern classes of stable Higgs extensions. The 
results in Section \ref{sect:Bog} generalize  the results 
described in \cite{DUW} for extensions of holomorphic bundles, 
with the proofs being one more illustration of how results for 
holomorphic bundles can be recast as results for Higgs bundles. 
Going one step further than in \cite{DUW}, we describe in detail 
the implications of attaining equality in the Bogomolov 
inequalities. 

\par Finally, in the Appendix, we extend to Higgs bundles 
the construction of Bott-Chern forms. These forms play an 
important role in the proof of the Hitchin-Kobayahi 
correspondence. In fact our proof uses only two special cases and 
all the requisite results can be extracted from the literature. 
The available treatments are however all somewhat ad hoc. We have 
thus undertaken a more systematic and general discussion, but have 
confined it to an Appendix.  Our results show how the original 
constructions of Bott and Chern for holomorphic bundles go over in 
their entirety to the case of Higgs bundles. This can be viewed as 
yet another illustration of the main underlying principle of this 
paper.  
\section{The Objects}\label{sectn:Objects}
Let $X$\ be a closed \kler\ manifold of dimension $d$ and with 
\kler\ form $\omega$. A Higgs sheaf (cf. \cite{S1, S2, S3, S4})
on $X$ is a pair $(\SE, \Theta)$ where $\SE$ is a coherent sheaf 
on $X$ and $\Theta$ is a morphism $\Theta: \SE \lra \SE \otimes 
\Omega^1_X$ (where $\Omega^1_X$ is the sheaf of holomorphic
sections of the cotangent bundle $T^*X$) 
such that 
$\Theta \wedge \Theta =0$. 
If $\cE$ is locally free, $\Theta$ can be thought of as a holomorphic
section of $\cE nd(\cE)\otimes \Omega^1_X$.
A morphism of Higgs sheaves $f: (\SE,\Theta) \lra 
(\SF,\Psi)$ is a morphism of sheaves $\overline{f}:\SE\lra 
\SF$ such that the following diagram commutes
\be
\begin{CD}\label{eq:HiggsMorf}
\SE @>{\Theta}>> \SE\otimes\Omega^1_X \\
@V{\overline{f}}VV           @V{\overline{f}\otimes id}VV \\ 
\SF @>{\Psi}>> \SF\otimes\Omega^1_X \\
\end{CD}
\end{equation}
Since the category of Higgs sheaves is abelian, the notion of 
exact sequence makes sense. 
\begin{definition}\label{defn: HiggsShExt}
An extension of Higgs sheaves (or Higgs extension) 
is a short exact sequence 
\be\label{eq:HiggsShExt}
\begin{CD}
0 @>>> (\cE_1,\Theta_1) @>{i}>> (\cE,\Theta) @>{q}>>
(\cE_2,\Theta_2) @>>> 0 \\
\end{CD}
\end{equation}
A morphism between extensions of Higgs sheaves is a commutative
diagram
\be\label{eq:HiggsShExtMorph}
\begin{CD}
0 @>>> (\cE'_1,\Theta'_1) @>>> (\cE',\Theta') @>>>
(\cE'_2,\Theta'_2) @>>> 0 \\
@. @VV{f_1}V  @VV{f}V  @VV{f_2}V @. \\
0 @>>> (\cE_1,\Theta_1) @>>> (\cE,\Theta) @>>>
(\cE_2,\Theta_2) @>>> 0 \\
\end{CD}
\end{equation}
\end{definition}
It follows that a morphism of Higgs extensions is
an isomorphism if and only if the three morphisms $f_1$, $f$
and $f_2$ are isomorphisms of Higgs bundles.
%
%
%
\section{Stability}\label{sectn:Stability}
The notions of stability for holomorphic bundles adapt 
straightforwardly to define both slope- and Gieseker stability for 
Higgs bundles (cf. \cite{S1,S2,S3,S4} and \cite{H}). In 
\cite{BGPcag} and \cite{DUW} these notions are defined for 
extensions of holomorphic bundles (or more generally, extensions 
of coherent sheaves). In this section we combine both of these to 
define stability for extensions of Higgs sheaves. As usual, the 
definition involves a numerical criterion on all subobjects. We 
must thus first define subobjects. 
\begin{definition}\label{def:subobj}
Consider a morphism of Higgs extensions
\be
\begin{CD}
0 @>>> (\cE'_1,\Theta'_1) @>>> (\cE',\Theta') @>>>
(\cE'_2,\Theta'_2) @>>> 0 \\
@. @VV{f_1}V  @VV{f}V  @VV{f_2}V @. \\
0 @>>> (\cE_1,\Theta_1) @>>> (\cE,\Theta) @>{q}>>
(\cE_2,\Theta_2) @>>> 0 \\
\end{CD}
\end{equation}
If $f_1$, $f$ and $f_2$ are injective, 
then the extension in the first row is called a subextension of the 
extension in the second row. A subextension is called 
proper if $\SE'$ is a proper subsheaf of $\SE$. 
\end{definition}
\smallskip
\remark 
Note that giving a proper subextension is the same thing as giving a 
proper subsheaf $\SE'$ of $\SE$  that is invariant under $\Theta$, 
in the sense that the image of 
$\Theta(\SE')$ is in $\SE'\otimes \Omega^1_X \subset \SE \otimes
\Omega^1_X$. Indeed, if $\SE'$ is invariant under $\Theta$, it
defines a Higgs subbundle $(\SE',\Theta')$, then 
we can recover $(\SE_2',\Theta'_2)$ as the
image of $\SE'$ under $q$, and $(\SE'_1,\Theta'_1)$ is recovered as
the kernel.

We can now define the notion of slope (or Mumford) stability. 
\begin{definition}[Slope stability]\label{def:aslope} 
Fix $\alpha<0$. Given a Higgs extension
\be
0 \lra (\SE_1,\Theta_1) \lra (\SE,\Theta) \lra (\SE_2,\Theta_2) 
\lra 0,
\end{equation}
define its $\alpha$-slope as 
\be\label{eq:aslope}
\mu_{\alpha}(\cE)=\mu(\cE)+\alpha\frac{\rk(\cE_2)}{\rk(\cE)}\ ,
\end{equation}
We say that a Higgs extension is $\alpha$-slope stable (resp. semistable), 
if for all proper 
subextensions, we have
\be\label{eq: astable}
\mu_{\alpha}(\cE')<\mu_{\alpha}(\cE)\ \qquad \text{(resp. $\leq$).}
\end{equation}
\end{definition}
\bigskip
\remark In particular, if $(\cE,\Theta)$\ is 
$\alpha$-stable then $\mu_{\alpha}(\cE_1)<\mu_{\alpha}(\cE)$.  
It follows from this that $\alpha>\mu(\cE_1)-\mu(\cE_2)$, i.e. the 
allowed range for the parameter $\alpha$ is 
\be\label{eq:arange}
\mu(\cE_1)-\mu(\cE_2)<\alpha<0\ .
\end{equation} 
In section \ref{sect:GIT}, where we construct moduli spaces, we 
will need a notion of Gieseker (semi)stability for Higgs extensions.
\begin{definition}[Gieseker stability]\label{defn:GiesStab}
Fix $\alpha<0$. 
Let $P(\SE,m)$\ denote the Hilbert 
polynomial of $\SE$\ .
A Higgs extension is called 
$\alpha$-Gieseker stable
(resp. semistable) if for all proper subextensions we have 
\begin{enumerate}
\item [(i)]
\be\label{eq: Gies1}
\mu_{\alpha}(\cE')\leq\mu_{\alpha}(\cE)\ .
\end{equation}
\item [(ii)]
If equality holds in (i), then 
\be\label{eq:Gies2}
\frac{P(\SE',m)}{\rk(\SE')} \leq \frac{P(\SE,m)}{\rk(\SE)}
\quad \text{for $m\gg 0$}
\end{equation}
\item [(iii)]
If equality holds in (i) and (ii), then 
\be\label{eq:Gies3}
\frac{P(\SE_2',m)}{\rk(\SE_2')} > \frac{P(\SE_2,m)}{\rk(\SE_2)}
\quad \text{(resp. $\geq$) for $m\gg 0$}
\end{equation}
\end{enumerate}
\end{definition}
As usual, we have the following implications 
$$
\text{$\alpha$-slope stable $\Longrightarrow$ $\alpha$-Gieseker stable
$\Longrightarrow$}
$$
$$
\text{$\Longrightarrow$ $\alpha$-Gieseker semistable $\Longrightarrow$ 
$\alpha$-slope semistable}
$$
\section{Differential Geometric Description and Metric Equations}
\label{sectn:Eqtns}
All the essential differential geometric machinery for Higgs 
bundles can be found in \cite{S3,S4} and \cite{H}. We thus give 
only a brief summary, emphasizing the aspects needed later in this 
paper.  Denoting the underlying smooth bundle of a holomorphic 
bundle $\cE$\ by $E$, we can describe the holomorphic structure on 
$\cE$ by an integrable partial connection, i.e. by a $\bC$-linear 
map 
\be\label{eq:dbareDefn}
\dbare:\Omega^0(E)\lra\Omega^{0,1}(E)
\end{equation}
\noi which satisfies the $\overline{\partial}$-Leibniz formula and 
also the integrability condition 
\be
\dbare\circ\dbare=\dbare^2=0
\end{equation} 
\noi A Higgs bundle $(\cE,\Theta)$ can thus be specified by a triple 
$(E,\dbare,\Theta)$\ where 
\begin{itemize}
\item $E$ is a smooth complex bundle on $X$,
\item $\dbare:\Omega^0(E)\lra\Omega^{0,1}(E)$ satisfies 
the $\overline{\partial}$-Leibniz formula and $\dbare^2=0$, 
\item $\Theta\in\Omega^{1,0}(End(E))$ satisfies 
$\dbare(\Theta)=0$\ and $\Theta\wedge\Theta=0$
\end{itemize}
\noi Instead of treating the holomorphic structure ($\dbare$) and 
the Higgs field ($\Theta$) as separate, we can combine them to 
define the Higgs operator 
\be\label{eq:nabla''}
 \nabla''=\dbare+\Theta:\Omega^0(E)\lra\Omega^{0,1}(E)\oplus\Omega^{1,0}(E)
\end{equation}
 
\noi Notice that this differs from the partial connection $\dbare$ in that its image 
is not confined to $\Omega^{0,1}(E)$. However, like  $\dbare$, it 
satisfies the $\dbar$-Leibniz formula and extends in the usual way 
to an operator on $\Omega^p(E)$.  Conversely, given any 
$\bC$-linear map 
$\nabla'':\Omega^0(E)\lra\Omega^{1}(E)$\ which satisfies the 
$\dbar$-Leibniz formula, we can separate it into  
$\nabla''=\dbare+\Theta$, corresponding to the splitting 
$\Omega(E)^1=\Omega^{0,1}(E)\oplus\Omega^{1,0}(E)$. 
\noi The  integrability condition, 
\be\label{eq:HiggsInt}
(\nabla'')^2=0\ , 
\end{equation}
\noi is clearly equivalent to the defining conditions of a Higgs bundle, 
viz. 
\bes
(\dbare)^2=0\ ,\ \dbare(\Theta)=0\ ,\ \Theta\wedge\Theta=0\ . 
\end{equation*}
\noi We thus arrive at the following description of a Higgs 
bundle, formally identical to the differential geometric 
description of a holomorphic bundle, but with the operator 
$\dbare$ replaced by the operator $\nabla''$. 
\begin{definition}[Higgs operator description]\label{defn:HiggsOpDesc}
A Higgs bundle on $X$is a pair $(E,\nabla'')$ in which $E$\ is a 
smooth bundle on $X$ and 
$\nabla'':\Omega^0(E)\lra\Omega^{1}(E)$ is a $\bC$-linear map
which satisfies the $\dbar$-Leibniz formula and the integrabiltiy 
condition (\ref{eq:HiggsInt}). 
\end{definition} 
Given a Hermitian bundle metric, $H$, on $E$, we can complete 
$\nabla''$\ so as to define a connection. To do so, we first define the adjoint 
$\Theta^*_H\in \Omega^{0,1}(EndE)$ by the condition 
that for all sections $s,t\in\Omega^0(E)$ 
\be\label{eq:Theta*def}
(\Theta s,t)_H=(s,\Theta^*_H t)_H\ . 
\end{equation}
\noi If we fix a local frame $\{e_i\}$\ for $E$, and define the Hermitian matrix
\be\label{eq:Hij}
H_{ji}=(e_i,e_j)_H\ , 
\end{equation}
\noi then $\Theta^*_H$\ is represented by the matrix
\be\label{eq:Theta*}
\Theta^*_H=H^{-1}\overline{\Theta}^TH\ .
\end{equation}
\noindent More explicitly, if we write
\be
\Theta=\sum_{\alpha}[\Theta^{\alpha}]_{ij} \otimes
\omega_{\alpha}\ ,
\end{equation}
\noindent where the $\omega_{\alpha}$\ are $(1,0)$-forms and the
matrices $[\Theta^{\alpha}]_{ij}$\ are local descriptions (with 
respect to the frame $\{e_i\}$) of bundle endomorphisms, then 
\be \Theta^*_H=\sum_{\alpha}[\Theta^{*,\alpha}_H]_{ij} \otimes
\overline{\omega}_{\alpha}\ ,
\end{equation}
\noindent where
\be\label{[Theta*]}
[\Theta^{*,\alpha}_H]_{ij}= 
H^{-1}_{ip}[\overline{\Theta^{*,\alpha}_H]}^T_{pq}H_{qj} \ . 
\end{equation} 
 
\begin{definition}\label{defn:HiggsConn}Define 
\be\label{eq:nabla'}
\nabla'_H=D'_H + \Theta^*_H\ .
\end{equation}
\noindent where $D(\dbare,H)=\dbare+D'_H$\ is the Chern connection
compatible with $\dbare$\ and $H$. The Higgs Connection is then 
defined by 
\be\label{eq:nabla}
\nabla =\nabla''+\nabla'_H \ .
\end{equation}
\noi The curvature of this connection 
\be\label{eq:HiggsCurv}
\Fhiggs = \nabla^2\ ,
\end{equation} 
\noi is called the Higgs curvature.
\end{definition}
\noindent \remark The Higgs curvature, like the curvature of any connection,
is a section of $\Omega^2(M,EndE)$. Unlike in the case of the 
Chern connection, $\Fhiggs$ does not have complex form type 
$(1,1)$.   The Higgs connection and its curvature do however have 
the following two crucial features: 
\begin{itemize}
\item (Kahler identities) 
\be\label{eq:KlerId}
i[\Lambda, \nabla'']=(\nabla'_H)^*\ ,\ i[\Lambda, 
\nabla'_H]=-(\nabla'')^*\ ,
\end{equation}
\noi where the adjoints are taken with respect to the metric $H$\ 
and
\be\label{eq:Lambda}
\Lambda:\Omega^{p,q}(E)\lra\Omega^{p-1,q-1}(E)
\end{equation}
is the adjoint of wedging with the \kler\ form on $X$.
\item (Bianchi identity) 
\be\label{eq:BianId}
\nabla'_H(\Fhiggs)=0=\nabla''(\Fhiggs)\ .
\end{equation}
\end{itemize}
\noindent Notice that these are direct analogs of the properties
enjoyed by the Chern connection, with $\nabla''$ and $\nabla'_H$\ 
playing the role here that $\dbare$\ and $D'_H$\ play for the 
Chern connection. This formal correspondence, which leads directly 
to the underlying principle mentioned in the Introduction, is 
summarized in Table 1.

\begin{table}    
\begin{center}
\begin{tabular}{|c|c|c|}
\hline
\hline
& &  \\
& Holomorphic bundle &\ Higgs bundle\\
& & \\
\hline
& &  \\
\text{underlying smooth bundle}&\ $E$ &\ $E$ \\
& & \\
\hline
& &  \\
\text{differential operator} & 
$\dbare:\Omega^0(E)\lra\Omega^{0,1}(E)$ &
$\nabla'':\Omega^0(E)\lra\Omega^{1}(E)$\\
& & \\
\hline
& &  \\
\text{integrability condition} & $\dbare^2=0$\ &\ $(\nabla'')^2=0$\\
& &  \\
\hline
& &  \\
\text{complementary operator}\ &\ $(D'_H)^*=i[\Lambda,\dbare]$\ &\ 
$(\nabla'_H)^*=i[\Lambda,\nabla'']$\\
& &  \\
\hline
& &  \\
\text{connection}\ &\ $D=\dbare+D'_H$\ &\ 
$\nabla=\nabla''+\nabla'_H$\\
& & \\
\hline
& &  \\
\text{ gauge theory equations}\ & & \\
\text{for special metrics}\ & $i\Lambda F^D_H=\mu\bI$\ & 
$i\Lambda \Fhiggs =\mu\bI$\\
& & \\
\hline
\hline
& &  \\
\text{(other) Kahler identity}\ &\ $(\dbare)^*=-i[\Lambda,D'_H]$\ &\ 
$(\nabla'')^*=-i[\Lambda,\nabla'_H]$\\ 
& &  \\
\hline
& &  \\
\text{Bianchi curvature identities}\ &\ $\dbare({F^D_H})=D'_H(F^D_H)=0$\ &\ 
$\nabla''(\Fhiggs)=\nabla'_H(\Fhiggs)=0$\\
& &  \\
\hline
\hline
\end{tabular}
\begin{caption}
{}Differential Geometric Dictionary, illustrating the formal 
similarity resulting from using the Higgs operator 
$\nabla''=\dbare+\Theta$ to encode the Higgs structure 
in a Higgs bundle 
\end{caption}
\end{center}
\end{table}

We now consider an extension of Higgs bundles,  
\bes
\Hextn
\end{equation*}
\noi i.e. a Higgs 
extension as in Definition \ref{defn: HiggsShExt} but in which the 
sheaves are locally free.  If we denote the underlying smooth 
bundle of 
$\cE$\ by $E$, then we can fix a smooth splitting 
$E=E_1\oplus E_2$, where the summands are the underlying smooth bundles 
for $\cE_1$\ and $\cE_2$. Thus the sub-Higgs bundle in the 
extension is described by the triple 
$(E_1,\delbar_1,\Theta_1)$, and the quotient Higgs bundle by 
$(E_2,\delbar_2,\Theta_2)$. The Higgs extension is 
then specified by the triple $(E,\dbare,\Theta)$\ where 
\begin{itemize}
\item the holomorphic structure is of the form
\be\label{eq:dbare}
\dbare=\left(\begin{array}{cc}
\dbar_1 & \beta\\ 0 &\dbar_2 
\end{array}\right)\ ,\ 
\beta \ \text{a holomorphic section in}\ 
\Omega^{0,1}(Hom(E_2,E_1))\ ,
\end{equation}
\item and the Higgs field is of the form
\be\label{eq:Theta}
\Theta=\left(\begin{array}{cc}
 \Theta_1 & b\\ 0 &\Theta_2 
\end{array}\right)\ , 
\ b\ \text{a holomorphic section in}\ 
\Omega^{1,0}(Hom(E_2,E_1))\ .
\end{equation}
\end{itemize}
\noi Here the holomorphic 
structure on $Hom(E_2,E_1)$ is that induced by 
$\overline{\partial}_1$\ and $\overline{\partial}_2$. Alternatively, 
using Higgs operators to describe the Higgs bundles, we have
\bes
\HopExtn
\end{equation*}
\noi where, with respect to a smooth splitting $E=E_1\oplus E_2$, 
the Higgs operator on $E$ is of the form 
\be\label{eq:extnabla''}
\nabla''=\left(\begin{array}{cc}
 \nabla_1'' & b+\beta\\ 0 & \nabla_2'' \\
\end{array}\right)\ 
\end{equation}
\noi Suppose now that we have a metric $H$\ on the middle bundle in the 
extension. It then makes sense to talk of an orthogonal splitting 
$E=E_1\oplus E_2$. We can thus define a bundle automorphism 
$\bT:E\lra E$ which, with respect to the $H$-orthogonal splitting, 
is given by the matrix
\be\label{eq:T}
\bT=\lb\begin{array}{cc}
 \frac{n_2}{n}\bI_1 & 0\\ 0 &
-\frac{n_1}{n}\bI_2 
\end{array}\rb\ .
\end{equation}
\noi Here $n=\rk(E)$\ and $n_i=\rk(E_i)$.  We can now formulate the following gauge 
theoretic equations:
\begin{definition}\label{defn:aHEeqtn} Fix the real number 
$\alpha$. We say the metric $H$\ satisfies the $\alpha$-Higgs-
Hermitian-Einstein (\aHHE) condition if
\be\label{eq:defHE}
i\Lambda \Fhiggs=\mu\bI+\alpha\bT\ , 
\end{equation}
\noi where $\Fhiggs$\ is the Higgs curvature as in (\ref{eq:HiggsCurv}), 
$\Lambda$ is as in (\ref{eq:Lambda}), $T$\ is the bundle automorphism defined in
(\ref{eq:T}) and $\mu=\mu(\cE)$\ is the slope of $\cE$.
\end{definition}
\noi\remarks \begin{itemize}
\item In the case $\Theta=0$, when $\nabla''=\dbare$\ and thus
the Higgs curvature $\Fhiggs$ reduces to $F^D_H$ (the curvature of 
the Chern connection compatible with $H$\ and $\dbare$ on 
$E$), equation (\ref{eq:defHE}) becomes the deformed 
Hermitian-Einstein equation defined 
in \cite{BGPcag} on extensions of holomorphic bundles. 
\item If we set $\alpha=0$ then we recover the usual Higgs equation 
(defined by Simpson and Hitchin) for a metric on the Higgs bundle 
$(\cE,\Theta)$
\item Using the fact that $(\nabla'')^2=0$, we can express $\Lambda\Fhiggs$
as 
\be
\Lambda \Fhiggs=\Lambda (F^D_H + [\Theta,\Theta^*])\ ,
\end{equation}
\noi where $F^D_H$ is the curvature of the Chern connection.  
The $\alpha$-Higgs-Hermitian-Einstein equation can thus also be 
written in the form
\be
i\Lambda ( F^D_H + [\Theta,\Theta^*])=\mu\bI+\alpha\bT\ .
\end{equation}
\end{itemize}
\section{The Hitchin-Kobayashi Correspondence}
\label{sectn:HKcorresp}
In this section we investigate the relation between the 
$\alpha$-stability of a Higgs extension and the existence of a metric
satisfying the \aHHE\ condition. As in \S 4, we fix an extension 
of Higgs bundles 
\be\label{eq:Hext}
\Hextn
\end{equation}
\noi The underlying smooth bundles are denoted, as usual, by 
$E_1,\  E_2$, and $E$. With Higgs operators defined as in (\ref{eq:nabla''})
we can thus equivalently describe the extension as
 
\be\label{eq:HopExt}
\HopExtn
\end{equation}
The Hitchin-Kobayashi correspondence asserts that 
$\alpha$-stability is equivalent to the existence of an \aHHE\ metric. 
In Section \ref{subs:easyHK} we prove that existence of an \aHHE\ 
metric implies $\alpha$-(poly)stability. The converse is proved in 
Section \ref{subs:hardHK}. In both cases we see the advantage of 
encoding the Higgs structure in the Higgs operator; having done 
so, the proofs amounts to little more than using the dictionary 
provided in Table 1 to adapt the corresponding proofs for 
extensions of holomorphic bundles (as in \cite{BGPcag}). 
\subsection{The Easy Direction}\label{subs:easyHK}
\bigskip
\begin{thm}\label{thm:easyHK} Fix $\alpha<0$.  Suppose that the Higgs extension 
(\ref{eq:Hext}) supports a metric with respect to which the smooth 
splitting $E=E_1\oplus E_2$ is orthogonal, and satisfying the 
\aHHE\ condition (\ref {eq:defHE}). Then either the Higgs 
extension is 
$\alpha$-stable or it splits as a direct sum of $\alpha$-stable 
Higgs extensions, all with the same $\alpha$-slope. 
 
\end{thm}
\pf Suppose that the metric $H=H_1\oplus H_2$\ on $E$\ 
satisfies (\ref{eq:defHE}).  Let $\nabla=\nabla''+\nabla'_H$\ be 
the Higgs connection determined by $H$\ and the Higgs operator on 
$E$, and let $\Fhiggs$ be its curvature (as in Definition 
\ref{defn:HiggsConn}).
Let $\cE'\subset\cE$\ be any Higgs subsheaf, with corresponding 
Higgs subextension 
\be
\Hsextn
\end{equation}
\noindent If $\cE'$\ is a saturated subsheaf then it is locally free outside 
of a codimension two subset, say $\Sigma$ in $X$. We can thus 
define a projection $\pi:\cE|_{X-\Sigma}\lra\cE'|_{X-\Sigma}$. 
Since 
$(\cE',\Theta')$\ is a Higgs subsheaf, we can compute the degree of 
$\cE'$ by the formula (cf. \cite{S3}, Lemma 3.2) 
\be\label{eq:degree}
deg(\cE')=i\int_XTr(\Lambda \pi 
\Fhiggs)-\int_X|\nabla''\pi|^2_H\ . 
\end{equation}
\noindent But by (\ref{eq:defHE})
\be
i\Lambda \Fhiggs=\lb\begin{array}{cc} 
\tau_1\bI_1 & 0\\ 0 &
\tau_2\bI_2 
\end{array}\rb\ ,
\end{equation}
\noindent where 
\be\label{eq:tau}
\begin{array}{cl}
\tau_1&=\mu+\alpha\frac{n_2}{n}\ ,\\ 
\tau_2&=\mu-\alpha\frac{n_1}{n}\ .
\end{array}
\end{equation} 
\noindent   It follows (precisely as in Proposition 3.8 of \cite{BGPcag}) that 
\be\label{eq:intF}
i\int_XTr(\Lambda \pi \Fhiggs)= n'_1\tau_1+n'_2\tau_2\ , 
\end{equation}
\noindent where $n'_1=rank(\cE'_1)$ and $n'_2=rank(\cE'_2)$. Notice that 
the first of the relations in (\ref{eq:tau}) can be written as 
$\tau_1=\mu_{\alpha}(\cE)$, and that together they imply 
$\alpha=\tau_1-\tau_2$. Combining (\ref{eq:intF}) and (\ref{eq:degree}) thus leads to 
\be
\mu_{\alpha}(\cE')=\mu_{\alpha}(\cE) -\int_X|\nabla''\pi|^2_H\ ,
\end{equation}
\noi from which the conclusion follows in the usual way.
\qed
\subsection{The Hard Direction}\label{subs:hardHK}
We now consider the converse of Theorem \ref{thm:easyHK}. Keeping 
the notation of Section \ref{subs:easyHK}, we show that if a Higgs 
extension ({\ref{eq:Hext}}) is 
$\alpha$-stable, then $\cE$\ admits a metric with 
respect to which the smooth splitting $E=E_1\oplus E_2$ is 
orthogonal and satisfying the \aHHE\ equation (\ref{eq:defHE}), 
i.e. such that 
\bes
i\Lambda \Fhiggs=\mu\bI+\alpha\bT\ . 
\end{equation*}
\noi As in \cite{S3} and \cite{BGPcag}, we can separate the trace and trace-free 
parts of this equation.  We can always fix $det(H)$\ so that 
\be\label{eq:TrHE}
i\Lambda Tr(\Fhiggs)=n\mu\ . 
\end{equation}
\noi In fact, since $[\Theta,\Theta^*]=0$ has zero trace, 
$i\Lambda Tr(\Fhiggs)$\ is the same for the Higgs connection as it is 
for the (metric) Chern connection. The above equation is thus 
satisfied if $det(H)$\ is the Hermitian-Einstein metric on the 
determinant line bundle $det(\cE)$.  Henceforth, we assume that we 
have fixed a background metric, $K$, such that $i\Lambda 
Tr(F_K)=n\mu$. 
It remains therefore to prove that $E$\ admits a metric satisfying  
\be\label{eq:TrfreeHE}
i\Lambda F^{\perp}_H=\alpha\bT\ , 
\end{equation}
\noi where $F^{\perp}=F-\frac{1}{n}Tr(F)\bI$\ is the trace-free
part of $F$. 
\bigskip
The proof follows the standard pattern for Hitchin-Kobayashi 
correspondences. The method we use is essentially that of Simpson, 
with modifications as in \cite{BGPcag} to accommodate the features 
arising from the extension structure (i.e. the non-zero right hand 
side in the equation). We thus give only a sketch of the proof, in 
which we fully describe all novel modifications, but do not repeat 
the details that can be found in \cite{BGPcag}, 
\cite{S3} and \cite{Do1}. 
Let 
\be\label{eq S(K)}
S(K)=\{s\in\Omega^0(X,EndE) | s^{*_K}=s\ ,\ Tr(s)=0\} \ . 
\end{equation}
\noi Then any other metric with the same determinant as $K$\ can be 
described by $Ke^s$, with $s\in S(K)$.  Fix an integer $p>2n$, and 
define 
\be\label{eq:Met}
\cM et^p_2=\{H=Ke^s \ | s\in L^p_2(S(K))\}\ .
\end{equation}
\noi We now define a Donaldson functional on $\cM et$\ whose critical 
points are solutions to (\ref{eq:TrfreeHE}).  The original 
Donaldson functional was defined using Bott-Chern forms for pairs 
of metrics, and had Hermitian-Einstein metrics on holomorphic 
bundles as its critical points. The generalization required to 
accommodate the extra structure of a Higgs bundle is due to 
Simpson, while the adaptation to the case of stable extensions can 
be found in 
\cite{BGPcag}]. Here we must combine both of these modifications. 
Given metrics 
$H$ and $K$, we denote the functional defined by Donaldson by
$M_D(K,H)$. It's definition in terms of Bott-Chern classes is
\be\label{eq:MD(H,K)}
M_D(H,K)=\int_XR_2(H,K)\wedge\omega^{d-1}\ , 
\end{equation}
\noi where $R_2$\ is the Bott-Chern form associated with the polynomial 
$-\frac{1}{2}Tr(AB+BA)$. Donaldson also gave a more 
explicit formula which applies for pairs $(H,K)$\ when 
$H=Ke^s$\ with $s\in S(K)$. Simpson's generalization of 
$M_D$\ can be obtained directly from this formula: one simply 
replaces the Chern connection by the Higgs connection. We will 
denote Simpson's functional by $M_S(H,K)$. Though it's not needed 
in this proof, and was not formulated in this way by Simpson, this 
modification can put in a more general framework. In the Appendix 
we show how it can be seen as the result of a modification of the 
Bott-Chern forms themselves.  
The functional used in \cite{BGPcag} for metrics on $E=E_1\oplus 
E_2$ can be defined as 
\be\label{eq:Mtau}
M_{\tau_1,\tau_2}(H,K)=M_D(H,K)-2(\tau_1-\tau_2)\int_X 
R_1(H_1,K_1) 
\wedge\omega^d\ ,
\end{equation}
\noi where $H_1$\ and $K_1$\ are the induced metrics on $E_1$ and the 
Bott-Chern form $R_1$ is given by 
\be\label{eq:R1}
R_1(H,K)=\log\det(K^{-1}H)=Tr(\log K^{-1}H)\ . 
\end{equation}
\noi We can combine this with Simpson's generalization if we 
replace $M_D$\ by $M_S$. We then get the following, which is the 
appropriate functional for extensions of Higgs bundles: 
\begin{definition}\label{def:MHiggs}
\be\label{eq:MHiggst}
M^{Higgs}_{\tau_1,\tau_2}(H,K)=M_S(H,K)-2(\tau_1-\tau_2)\int_X 
R_1(H_1,K_1) \wedge\omega^d\ , 
\end{equation}
\noi or, setting $\alpha=\tau_1-\tau_2$,
\be\label{eq:MHiggsa}
M^{Higgs}_{\alpha}(H,K)=M_S(H,K)-2\alpha\int_X R_1(H_1,K_1) 
\wedge\omega^d\ .
\end{equation}
\noi If we fix one of the metrics, say $K$ , we can define  
\be\label{eq:MHiggsa(H)}
M^{Higgs}_{\alpha}(H)=M^{Higgs}_{\alpha}(H,K). 
\end{equation}
\noi Following \cite{BGPcag}, we now define 
$m^0_{\alpha}:\cM et\longrightarrow \Omega^0(X,EndE)$\ by 
\be\label{eq:m0}
m^0_{\alpha}(H)=\Lambda F^{\perp}_H + i\alpha \bT_H\ , 
\end{equation}
\noi where, with respect to the $H$-orthogonal splitting 
$E=E_1\oplus E_2$\ ,
\be\label{eq:T2}
\bT_H= \lb\begin{array}{cc}
 \frac{n_2}{n}\bI_1& 0\\
0 &-\frac{n_1}{n}\bI_2 
\end{array}\rb
\end{equation}
\end{definition}
\noi The crucial properties of $M^{Higgs}_{\alpha}$\ and $m^0_{\alpha}$\ are described in 
the next proposition. 
\begin{prop}\label{prop:properties}
\begin{enumerate}
\item Given any three metrics $H,K,J$, we have 
\be\label{eq:additive}
M^{Higgs}_{\alpha}(H,K) 
+M^{Higgs}_{\alpha}(K,J)=M^{Higgs}_{\alpha}(H,J)\ . 
\end{equation}
\item  If $H(t)=He^{ts}$\ with $s\in S(H)$, then 
\be\label{eq:Mdot}
\frac{d}{dt} M^{Higgs}_{\alpha}(H(t))=
2i \int_X Tr\left(s m^0_{\alpha}(H(t))\right)\ . 
\end{equation} 
\item Define the operator $L$\ on $L^p_2(S(H))$\ by 
\be\label{eq:L}
L(s)=\frac{d}{dt}m^0_{\alpha}(H(t))|_{t=0}\ . 
\end{equation}
\noi If $s\in S(H)$\ is given by $s= 
\lb\begin{array}{cc}
s_1& u\\u^*&s_2 
\end{array}\rb$
with respect to the $H$-orthogonal splitting $E=E_1\oplus E_2$\ , 
and $H(t)=He^{ts}$, then 
\be\label{eq:Mddot}
\begin{array}{cl}
2i\langle s,L(s)\rangle_H&=\frac{d^2}{dt^2} 
M^{Higgs}_{\alpha}(H(t))|_{t=0}\\ 
&=\parallel \nabla''(s)\parallel^2_H-\alpha\parallel 
u\parallel^2_H\end{array} 
\end{equation}
\item If $s\in S(H)$\ and $K=He^s$, then 
\be\label{eg:Delta|s|}
\Delta |s|\le 2(|m^0_{\alpha}(H)|_H+|m^0_{\alpha}(K)|_K)\ ,
\end{equation}
\noi where the norm on $|s|$ can be with respect to either $H$\ or $K$. 
\end{enumerate}
\end{prop}
{\em Proof of 1. and 2.} When $\alpha=0$, these results follow as 
in \cite{S3} (\S 5) and \cite{Do2} (or, equivalently, follow from 
the properties of Bott-Chern forms, as described in the Appendix).  
The modification required when 
$\alpha<0$\ is exactly the same as described in the proof of 
Proposition 3.11 in \cite{BGPcag}. 
{\em Proof of 3.} The proof is formally identical to that in 
Proposition 3.11 in \cite{BGPcag}, except we replace the result 
about the second variation of $M_D$ with the corresponding result 
for 
$M_S$, viz. 
\be\label{eq:MSddot}
\frac{d^2}{dt^2}M_S(H(t))|_{t=0}=\parallel \nabla''(s)\parallel_H^2\ .
\end{equation}
\noi This result can be found in \cite{S3}. It can also be 
derived directly from the properties of Bott-Chern forms, as in 
Proposition \ref{prop: MHiggsderivs} of the Appendix. 
\bigskip
{\em Proof of 4.} When $\alpha=0$, this is part (d) of Lemma 3.1 
in \cite{S3}. In general we have 
\be\label{eq:deltam0}
m^0_{\alpha}(H)-m^0_{\alpha}(K)=(m^0_{0}(H)-m^0_0(K)) 
+i\alpha(\bT_H-\bT_K)\ . 
\end{equation} 
\noi This changes the computation in Simpson's proof by the introduction 
of an extra term of the form 
\be\label{eq:extraterm}
\alpha Tr\left(e^s(\bT_H-\bT_K)\right)\ .
\end{equation}
\noi But $Tr (e^s\bT_H)=Tr (e^s\bT_K)$, so the extra term does not 
affect the result. 
\qed

\begin{cor}\label{cor:ker}
Suppose that $\alpha<0$\ and (\ref{eq:Hext}) is an 
$\alpha$-stable extension.  Then 
\be\label{eq:ker}
Ker(L)=0\ , 
\end{equation}
\noi where $L$\ is the operator defined above on $L^p_2(S(H))$.
\end{cor}
\pf Suppose that $L(s)=0$ for some $s\ne 0$. Then 
by (\ref{eq:Mddot}) we have $\nabla''(s)=0=u$. Recall that with 
respect to the $H$-orthogonal splitting $E=E_1\oplus E_2$, the 
holomorphic structure and Higgs field on $E$\ are given by 
(\ref{eq:dbare}) and (\ref{eq:Theta}). Thus
\be\label{eq:nabla2}
\nabla''=\left(\begin{array}{cc}
\nabla''_1 & \beta+b\\ 0 &\nabla''_2 
\end{array}\right)\
\end{equation}
\noi Writing $s= 
\lb\begin{array}{cc}
s_1& u\\u^*&s_2 
\end{array}\rb$, where $s_i\in L^p_2(S(K_i))$\ and 
$u\in \Omega^0(X,Hom(E_2,E_1))$, 
we thus have $\nabla''_1(s_1)=\nabla''_2(s_2)=0$. But 
$\nabla''_i(s_i)=0$\ is equivalent to 
\be\label{eq:nablasi}
\dbar_i(s_i)=0\ \textrm{and}\ [\Theta_i,s_i]=0\ 
\end{equation}
\noi The eiegenspaces of 
$s$\ thus split the extension (\ref{eq:Hext}) into a direct sum of Higgs extensions. 
Since $Tr(s)=0$ there must be at least two such summands. But this 
violates the stability criterion, since the 
$\alpha$-slope inequality cannot be satisfied by both summands.  
\qed

\remark This same computation shows that for any path 
$H(t)=He^{ts}$\ with $s\in S(H)$, we get
\be\label{eq:convex}
\frac{d^2}{dt^2} M^{Higgs}_{\alpha}(H(t))>0\ ,
\end{equation}
\noi i.e. $M^{Higgs}_{\alpha}$\ is a convex functional. 
\bigskip
Next, we fix a positive real number $B$\ such that 
$\parallel m^0_{\alpha}(K)\parallel^p_{L^p}\le B$, where
\be
\parallel m^0_{\alpha}(K)\parallel^p_{L^p}=\int_{X}{|m^0_{\alpha}(K)|^p_K
dvol} 
\end{equation}
\noi and define 
\be\label{eq:MetB}
\cM et^p_2(B)=\{H\in\cM et^p_2 |\
\parallel m^0_{\alpha}(H)\parallel^p_{L^p}\le B\ \}\ .
\end{equation}
  
\begin{lemma}\label{lem:extrema}
 If the extension (\ref{eq:Hext}) is 
$\alpha$-stable, then there are no extrema of $M^{Higgs}_{\alpha}$
on the boundary of this constrained space, and  the minima occur 
at solutions to the metric equation 
$m^0_{\alpha}(H)=0$. 
\end{lemma}
\pf
(as in \cite{B1}, Lemma 3.4.2), in which Ker(L)=0 is the key) 
\qed

\noi We thus look for minima of $M^{Higgs}_{\alpha}(H)$\ on $\cM et^p_2(B)$. 
To show that minima do occur, we need 
\begin{prop}\label{prop:mainest}(3.14 in \cite{BGPcag}) 
Either (\ref{eq:Hext}) is not $\alpha$-stable or we can find 
positive constants 
$C_1$\ and $C_2$\ such that
\be\label{eq:mainest}
sup|s|<C_1 M^{Higgs}_{\alpha}(Ke^s)+C_2\ 
\end{equation}
\noi for all $Ke^s\in \cM et^p_2(B)$. 
\end{prop}
\bigskip
\remark\  This proposition describes what might be called 
the Donaldson-Uhlenbeck-Simpson-Yau (DUSY) Alternative: either one 
can produce a minimizing sequence for the functional $ 
M^{Higgs}_{\alpha}$ - and hence a solution to the metric equation 
-  or one can use the functional to produce a sequence which in 
the limit destabilizes the extension (\ref{eq:Hext}). 
\bigskip

{\em Sketch of Proof} One first shows that for metrics in the 
constrained set $\Met^p_2(B)$, the 
$C^0$ estimate given above is equivalent to a $C^1$\ estimate of the same  
type.  The proof of this uses (\ref{eg:Delta|s|}) in Proposition 
\ref{prop:properties}, but is otherwise identical to that in \cite{S3}
or \cite{B1}. One then supposes that no such $C^1$ estimate holds. 
It follows that one may find an unbounded sequence of constants 
$C_i$ and metrics $Ke^s_i\in \cM et^p_2(B)$ such that the estimate 
is violated. After normalizing the $s_i$, this produces a sequence 
$\{u_i\}\subset L^P_2(S(K))$\ such that 
$\parallel u_i\parallel_{L^1}=1$. This has a weakly convergent 
subsequence in  $L^2_1(S(K))$, with non-trivial limit denoted by 
$u_{\infty}$. One then shows that the eigenvalues of $u_{\infty}$ are constant
almost everywhere.  This is done, as in (\cite{S3} \S 5), by 
making use of an estimate of the form: 
\begin{lemma}\label{lemma:3.13BGP} [Lemma 3.13, \cite{BGPcag}]
Suppose that $\alpha<0$\ and let $H=Ke^s$\ with $s\in 
L^p_2(S(K))$. Let $s= 
\lb\begin{array}{cc}
 s_1& u\\u^*&s_2 
\end{array}\rb$\ be the block decomposition of $s$\ with 
respect to the $K$-orthogonal splitting $E=E_1\oplus E_2$\ . Let 
$\Psi:\bR\times\bR\longrightarrow \bR$\ be the smooth 
function as in \cite{B1} (or \cite{S3}). Then 
\be\label{eq:Mest}
\begin{array}{cl}
M^{Higgs}_{\alpha}(H)&= i\int_X{Tr(s\Lambda F_K)}+ 
\int_x{(\Psi(s)\nabla'' s,\nabla'' s)_K}
-2\alpha R_1(H_1,K_1)\\ 
&\ge\ i\int_X{Tr(s\Lambda F_K)}+
\int_x{(\Psi(s)\nabla'' s,\nabla'' s)_K}-\alpha\int_x{Tr(s_1)}
\end{array}
\end{equation}
\noi where the meaning of $\Psi(s)$\ is as in \cite{B1}\ or \cite{S3}. 
\end{lemma}
\pf As in \cite{BGPcag}:  The first line follows from the computations in 
\cite{S3}.  The second uses the convexity properties of the function 
$R_1(H(t)_1,K_1)$, and  the fact that its first derivative at 
$t=0$\ is given by $\int_X Tr(s_1)$.
\qed

\bigskip
\noi Following the analysis in \cite{S3} ( Lemma 5.4), this leads to
\begin{prop}\label{prop:3.15BGP}(3.15 in \cite{BGPcag}) 
Let $\cF:\bR\times \bR\longrightarrow \bR$\ be  any smooth 
positive function which satisfies $\cF(x,y)\le 1/(x-y)$\ whenever 
$x>y$.  Then
\be\label{eq:uinftyest}
 i\int_X{Tr(u_{\infty}\Lambda F_K)}+
\int_x{(\cF(u_{\infty})\nabla'' u_{\infty},\nabla''
u_{\infty})_K}-\alpha\int_x{Tr(u_{\infty,1})}\le 0\ , 
\end{equation}
 
\noi where $u_{\infty}=\lb\begin{array}{cc}
 u_{\infty,1}& {*} \\ {*} & {*}
\end{array}\rb$\ with 
respect to the $K$-orthogonal splitting of $E$. 
\end{prop}
Since $Tr(u_{\infty})=0$, there are at least two distinct 
eigenvalues. Let 
$\lambda_1<\lambda_2,\dots,<\lambda_k$\ denote the distinct eigenvalues.
Setting $a_i=\lambda_{i+1}-\lambda_i$, one can thus define 
projections 
$\pi_i\in L^2_1(S(K))$\ such that
\be\label{eq:pii}
u_{\infty}=\lambda_r\bI-\sum_{i}^{k-1}{a_i\pi_i} 
\end{equation}
\begin{lemma}\label{lemma:pii}
The projections $\pi_i$\ satisfy 
\begin{enumerate}
\item $\pi_i\in L^2_1(S(K))$,
\item $\pi_i^2=\pi_i$
\item $(1-\pi_1)\nabla''(\pi_i)=0$
\end{enumerate}
\end{lemma}
\pf The $\alpha=0$\ case is proved in \cite{S3} (Lemma 5.6 
and succeeding remarks). The presence of the extra term depending 
on $\alpha$ in (\ref{eq:uinftyest}) does not affect the method of 
proof. 
\qed

\noi Each $\pi_i$\ thus defines a weak Higgs subbundle in the sense of 
Uhlenbeck and Yau \cite{UY}, as adapted by Simpson (\cite{S3}) for 
Higgs bundles, and hence produces a filtration of $\cE$\ by 
reflexive Higgs subsheaves 
\be
\cE_1\subset\cE_2\subset\dots\subset\cE_k=\cE
\end{equation}
\noi Each Higgs subsheaf $\cE_j$\ determines a Higgs subextension 
\be
0\longrightarrow \cE_{1,j}\longrightarrow \cE_{j} 
\longrightarrow \cE_{2,j}\longrightarrow  0\ .
\end{equation}
\noi Now define the numerical quantity 
\be\label{eq:Q}
Q=\lambda_k(r\mu(\cE)-r_1\tau_1-r_2\tau_2)- 
\sum_{i}^{k_1}{a_i(r_i\mu(\cE_i)-
r_{1,i}\tau_1-r_{2,i}\tau_2)}\ , 
\end{equation}
\noi  where
$\mu(\cE_i)$\ is the slope of $\cE_j$, and 
$r_{a,i}$\ is the rank of $\cE_{a,i}$.  Using Lemma \ref{lemma:3.13BGP} 
and the fact that 
$u_{\infty}=\lambda_r\bI-\sum_{i}^{k-1}{a_i\pi_i}$, 
one shows (by precisely the method in \cite{S3}) that $Q\le 0$. On 
the other hand, $\tau_1$\ and $\tau_2$\ are related by 
$r\mu(\cE)-r_1\tau_1-r_2\tau_2=0$, and if $(\ref{eq:Hext})$\ is 
$\alpha$-stable, then
\be\label{eq:rimu}
r_i\mu(\cE_i)-r_{1,i}\tau_1-r_{2,i}\tau_2<0 
\end{equation}
\noi for all $i=1,\dots, k-1$. Thus $Q$\ must be strictly positive if (\ref{eq:Hext}) is 
$\alpha$-stable. We conclude therefore that if (\ref{eq:Hext}) is 
$\alpha$-stable then there must be constants $C_1$ and $C_2$\ such 
that the estimate (\ref{eq:mainest}) holds. 
\qed

\bigskip
We can now prove 
\begin{thm}\label{thm:HardHK}
Fix $\alpha<0$\ and suppose that the Higgs extension 
(\ref{eq:Hext}) is $\alpha$-stable. Then $E$\ admits a unique 
metric $H$ with respect to which the smooth splitting $E=E_1\oplus 
E_2$ is orthogonal, with $det(H)=det(K)$, and such that 
\be
i\Lambda F^{\perp}_H=\alpha\bT\ . 
\end{equation}
\end{thm}
\pf By Proposition \ref{prop:mainest}, there is an estimate of the 
form in (\ref{eq:mainest}) and hence the functional 
$M^{Higgs}_{\alpha}$\ is bounded below. By Lemma 
\ref{lem:extrema}, a minimizing sequence produces a solution in 
$Met^p_2(B)$ to the equation 
$m^0_{\alpha}(H)=0$. The smoothness and uniqueness
of the solution follows in exactly the same way as in \cite{Do1},
\cite{S3} or \cite{B1}.  The smoothness is a result of elliptic regularity,
while the uniqueness is a consequence of the convexity properties 
of $M^{Higgs}_{\alpha}$. 
\qed
\section{Bogomolov Inequality}\label{sect:Bog}
The existence of a solution to the 
$\alpha$-Higgs-Hermitian-Einstein equations on an $\alpha$-stable
Higgs extension can be used to deduce topological constraints.  
The constraints are expressed as inequalities involving the Chern 
classes of the underlying bundles. As such, they are direct 
generalizations of the Bogomolov inequalities for stable 
holomorphic bundles. The notation in this section is as follows: 
\begin{itemize}
\item As in Section \ref{sectn:HKcorresp}, $(E,\nabla'')$ is a Higgs bundle 
which has the structure of an extension of Higgs bundles as in 
(\ref{eq:HopExt}), i.e. which can be written as 
\bes
\HopExtn\ .
\end{equation*}
\item The ranks of the underlying smooth bundles $E_1,E_2$\ and $E$\ 
are denoted by 
$n_1,n_2$\ and $n$\ respectively. 
\item The base space is the \kler\ manifold 
$(X,\omega)$. The dimension of $X$ is $d$, and its volume is $V$. 
\item Using the \kler\ form $\omega$\ and the chern classes $c_1(E),\ c_2(E)$,
we define the following characteristic numbers 
\be\label{CharNo}
C_2(E)=\int_Xc_2(E)\wedge\omega^{d-2}\quad\ ,\ 
C^2_1(E)=\int_Xc^2_1(E)\wedge\omega^{d-2}
\end{equation}
\end{itemize}
\noi With this notation, we prove the following results:
\begin{thm}\label{thm:BogIneq}(Bogomolov Inequality) Let $(E,\nabla'')$ be
a Higgs bundle which has the structure of an extension of Higgs 
bundles as in (\ref{eq:Hext}), i.e. which can be written as 
\bes
\HopExtn\ .
\end{equation*}
\noi Suppose that $(E,\nabla'')$\ is $\alpha$-polystable as 
an extension of Higgs bundles, for some $\alpha<0$. Then
\be\label{eq:BogIneq}
2C_2(E)-\frac{n-1}{n}C_1^2(E)+\alpha^2(\frac{n_1n_2}{n}) 
\frac{V(d-1)!}{4\pi^2d}\ge 0\ .
\end{equation}
\end{thm}
\bigskip
\begin{thm}\label{thm:BogEq} Let $(E,\nabla'')$ be as in 
Theorem \ref{thm:BogIneq}. Suppose that
$(E,\nabla'')$\ is $\alpha$-polystable as an extension of Higgs bundles
and that equality holds in (\ref{eq:BogIneq}), i.e.  its Chern 
classes satisfy 
\be\label{eq:BogEq}
2C_2(E)-\frac{n-1}{n}C_1^2(E)+\alpha^2(\frac{n_1n_2}{n})\frac 
{V(d-1)!}{4\pi^2d}= 0\ 
\end{equation}
\noi Then
\begin{enumerate}
\item with respect to the splitting $E=E_1\oplus E_2$\ we 
have 
\be\label{eq:dbare2}
\nabla''=\left(\begin{array}{cc} 
\nabla''_1 & 0\\ 0 & \nabla''_2
\end{array}\right)\ ,\ \mathrm{i.e.}\ 
\dbare=\left(\begin{array}{cc}
\dbar_1 & 0\\ 0 & \dbar_2
\end{array}\right)\ \mathrm{and}\
\Theta=\left(\begin{array}{cc}
\Theta_1 & 0\\ 0 & \Theta_2\ 
\end{array}\right)\ ,
\end{equation}
\item there is a metric $H=H_1\oplus H_2$, such that each 
summand satisfies 
\be\label{eq:Fperp}
F^{\perp}_{H_i}=0\ , 
\end{equation}
\noi and
\be\label{eq:TraceEq}
\frac{Tr(F_{H_1})}{n_1} - \frac{Tr(F_{H_2})}{n_2}=
\Lambda(\frac{Tr(F_{H_1})}{n_1}-\frac{
Tr(F_{H_2})}{n_2})\frac{\omega}{d} 
\,  
\end{equation}
\item the parameter $\alpha$\ has the value
\be\label{alpha}
\alpha =\mu_1-\mu_2\ ,
\end{equation}
\noi where 
\be\label{mu(i)}
\mu_i=\frac{2\pi}{n_i}\int_X \Lambda c_1(E_i)
\frac{\omega^d}{d!}\ .
\end{equation}
\end{enumerate}
\bigskip
\noi Conversely, if conditions (1)-(3) apply, then the Higgs 
extension is $\alpha$-polystable and its chern classes satisfy the 
equality \ref{eq:BogEq}. 
\end{thm}
\remark  Conditions (1) and (2) in Theorem \ref{thm:BogEq} together imply
that $(E,\nabla'')$\ splits as a direct sum of polystable Higgs 
bundles. 
\bigskip
\noindent We require the following key technical result :
\begin{prop}\label{prop:Simpson}(\cite{S3}, \S 3) If $\Fhiggs$\ is the curvature of
the Higgs connection determined by metric $H$\ on 
$(E,\nabla'')$, then
\be\label{eq:TrFF}
Tr(\Fhiggs\wedge \Fhiggs\wedge\omega^{d-2})= 
|\Fhiggs-\frac{1}{d}(\Lambda 
\Fhiggs)\omega|^2\frac{\omega^d}{d(d-1)}-|\Lambda 
\Fhiggs|^2\frac{\omega^d}{d^2}\ 
\end{equation}
\noi where $d=dim(X)$.  Similarly, if $(\Fhiggs)^{\perp}=
\Fhiggs-\frac{1}{d}Tr(\Fhiggs)\bI$, then 
\be\label{eq:TrFFperp}
(F^{\perp}_H\wedge F^{\perp}_H\wedge\omega^{d-2})= 
|F^{\perp}_H-\frac{1}{n}(\Lambda 
F^{\perp}_H)\omega|^2\frac{\omega^d}{d(d-1)}-|\Lambda 
F^{\perp}_H|^2\frac{\omega^d}{d^2}\ 
\end{equation}
\end{prop}
\pf This uses the following features of Higgs
connections: 
\be\label{eq:F(1,1)}
(\Fhiggs)^{1,1}+((\Fhiggs)^{1,1})^{*_H}=0 
\end{equation}
\be\label{eq:F(2,0)}
(\Fhiggs)^{2,0}=((\Fhiggs)^{0,2})^{*_H} 
\end{equation}
\qed

\bigskip
\noi \textit{Proof of Theorem \ref{thm:BogIneq}}  If $(E,\nabla'')$\ 
is $\alpha$-ploystable, then (by Theorem \ref{thm:HardHK}) it has 
a metric satisfying the \aHHE\ equation (\ref{eq:defHE}). Taking 
the trace-free part, i.e. (\ref{eq:TrfreeHE}) , we get 
\begin{align}\label{eq:normFperp}
||\Lambda F^{\perp}_H||^2&=\int_X|\Lambda F^{\perp}_H|^2\ 
\frac{\omega^d}{d!}\\ 
&=\int_X|\alpha\bT |^2\
\frac{\omega^d}{d!}\\ 
&=\alpha^2\frac{n_1n_2}{n}V\ ,
\end{align}
\noi where $V$\ is the volume of $X$. Also, using the Chern-Weil formulae
for $ch_2(E)$\ and $c_1(E)$, plus the identity 
$ch_2=\frac{1}{2}c_1^2-c_2$, we get 
\be\label{eq:normFperp2}
\begin{array}{cl}
\frac{1}{4\pi^2}\int_XTr(F^{\perp}_H\wedge
F^{\perp}_H\wedge\omega^{d-2})&=\frac{1}{4\pi^2}\int_X(Tr(\Fhiggs\wedge 
\Fhiggs)-\frac{1}{n}Tr(\Fhiggs)\wedge 
Tr(\Fhiggs))\wedge\omega^{d-2}\\ 
&=\int_X (-2ch_2(E)+\frac{1}{n}c_1^2(E))\wedge\omega^{d-2}\\
&=\int_X (2c_2(E)-\frac{n-1}{n}c_1^2(E))\wedge\omega^{d-2}\ .
\end{array}
\end{equation}
\noi Equation (\ref{eq:TrFFperp}) thus yields
\be\label{eq:ChernNoEq}
2C_2(E)-\frac{n-1}{n}C_1^2(E)+\alpha^2(\frac{n_1n_2}{n})\frac{V(d-1)!}{4\pi^2d} 
=\frac{(d-2)!}{4\pi^2}||F^{\perp}_H-\frac{1}{d}(\Lambda 
F^{\perp}_H)\omega||^2\ , 
\end{equation}
\noi where $C_2(E)$\ and
$C^2_1(E)$\ are as in (\ref{CharNo}).  Theorem \ref{thm:BogIneq} follows 
directly from this. 
\qed

\bigskip
\textit{Proof of Theorem \ref{thm:BogEq}} Suppose that
$(E,\nabla'')$\ is $\alpha$-polystable as an extension of Higgs 
bundles, and that (\ref{eq:BogEq}) holds. As in the previous 
proof, we  may thus assume that $E$\ supports a metric 
$H=H_1\oplus H_2$\ which satisfies the trace-free \aHHE\ equation 
(\ref{eq:TrfreeHE}). It then follows from (\ref{eq:ChernNoEq}) 
that the trace free part of the curvature, i.e. 
$F^{\perp}_H$, satisfies 
\be\label{eq:Fperp2}
F^{\perp}_H=-i\alpha\bT \frac{\omega}{d}\ . 
\end{equation}
\noi Applying the Bianchi identity, viz. $\nabla(\Fhiggs)=0$, and the
fact that (cf. Lemma \ref{lemma:dphi}) 
$dTr(\Fhiggs)=Tr\nabla(\Fhiggs)$, we  get 
\be\label{eq:nablaT}
\nabla(T)=0\ .
\end{equation}
\noi It follows from this that the subbundles corresponding to
eigenvalues $\frac{n_2}{n}$\ and $-\frac{n_1}{n}$\ of $\bT$\ both 
give rise to Higgs subbundles of $(E,\nabla'')$. Alternatively, 
one can compute the covariant derivative 
$\nabla(T)$\ and observe directly from (\ref{eq:nablaT}) that 
$\nabla''$ (and hence $\dbare$\ and $\Theta$) must be as in (\ref{eq:dbare2}).  
Either way, we have 
\be\label{eq:Fsplit}
\Fhiggs=\left(\begin{array}{cc} F_{H_1} & 0\\ 0 & F_{H_2} 
\end{array}\right)
\end{equation}
\noi and hence
\be\label{eq:Fperpsplit}
F^{\perp}_H= \left(\begin{array}{cc} 
 F^{\perp}_{H_1} & 0\\ 0 & F^{\perp}_{H_2} 
\end{array}\right) + (\frac{Tr(F_1)}{n_1}-\frac{Tr(F_2)}{n_2})\bT\ , 
\end{equation}
\noi where $F^{\perp}_{H_i}=F_{H_i}-\frac{Tr(F_{H_i})}{n_i}$\ for
$i=1,2$.  Combining this with (\ref{eq:Fperp2}), we see that
\be\label{eq:Fperpsplit2} 
\left(\begin{array}{cc}
F^{\perp}_{H_1} & 0\\ 0 & F^{\perp}_{H_2} 
\end{array}\right) =(\ \frac{Tr(F_2)}{n_2}-\frac{Tr(F_1)}{n_1}-i\alpha
\frac{\omega}{d})\ \bT\ ,
\end{equation}
\noi i.e.
\be\label{eq:Fperp12}
\begin{array}{c}
F^{\perp}_{H_1}=\frac{n_2}{n}(\ 
\frac{Tr(F_2)}{n_2}-\frac{Tr(F_1)}{n_1}-i\alpha 
\frac{\omega}{d})\ \bI_1\ ,\\
F^{\perp}_{H_2}=-\frac{n_1}{n}(\ 
\frac{Tr(F_2)}{n_2}-\frac{Tr(F_1)}{n_1}-i\alpha
\frac{\omega}{d})\ \bI_2\ 
\end{array}
\end{equation}
\noi Taking the trace of either of these equations yields (\ref{eq:TraceEq}).
Then integrating over $X$ yields (\ref{alpha}). 
Conversely, suppose that (1)\ -\ (3) apply. Then 
(\ref{eq:Fperpsplit}) implies 
\be
F^{\perp}_H=-i\alpha\bT\frac{\omega}{d}=\Lambda 
F^{\perp}_H\frac{\omega}{d}\ , 
\end{equation} 
\noi and hence that the right hand side of (\ref{eq:ChernNoEq}) vanishes. Thus, with
$H=H_1\oplus H_2$, we see that
$i\Lambda \Fhiggs=\mu\bI+\alpha\bT$, as required.  It remains 
to verify (\ref{eq:BogEq}). 
We write, for $i=1,2$ 
\be\label{eq:c1}
 c_1(E_i)=\delta_i\omega+\beta_i\ ,
\end{equation}
\be\label{eq:c2}
c_2(E_i) = a_i\omega^2 + b_i\wedge\omega + c_i 
\end{equation}
\noi where $\delta_i, a_i\in\bR$\ and 
$\beta_i, b_i\in\Omega^{(1,1)}(X,\bR)$\ are
primitive forms, and $c_i\wedge\omega^{(d-2)}=0$.  The condition 
in (\ref{eq:TraceEq}) then becomes 
\be\label{eq:betaEq}
\frac{\beta_1}{n_1}-\frac{\beta_2}{n_2}=0\ .
\end{equation}
\noi Using the identities,
\be
c_2(E_1\oplus E_2)=c_2(E_1)+c_2(E_2)+c_1(E_1)\wedge c_1(E_2)\ , 
\end{equation}
\noi and
\be
c_1(E_1\oplus E_2)=c_1(E_1)+c_1(E_2)\ , 
\end{equation}
\noi we thus compute
\be\label{eq:ChernEq2}
\begin{array}{cl}
2C_2(E)-\frac{(n-1)}{n}C_1^2(E)&= (2(a_1+a_2+\delta_1\delta_2) 
-\frac{n-1}{n}(\delta_1+\delta_2)^2)\omega^d\\ 
&\quad +2(\beta_1\wedge\beta_2-\frac{n-1}{n}(\beta_1+\beta_2)^2)\wedge\omega^{(d-2)}\\
&=\sum_{i=1,2}(2C_2(E_i)-\frac{n_i-1}{n_i}C_1^2(E_i))\\
&\quad\quad +\frac{n_1n_2}{n}(\frac{\delta_1^2}{n_1^2}+\frac{\delta_2^2}{n_2^2}-
2\frac{\delta_1\delta_2}{n_1n_2}) 
-\frac{n_1n_2}{n}(\frac{\beta_1}{n_1}-\frac{\beta_2}{n_2})^2 
\end{array}
\end{equation}
\noi By the Bogomolov inequality for polystable bundles, we have
\be\label{eq:bndlBog}
2C_2(E_i)-\frac{n_i-1}{n_i}C_1^2(E_i)= 0\ . 
\end{equation}
\noi Together with (\ref{eq:betaEq}), equation (\ref{eq:ChernEq2}) thus reduces to
\be
\begin{array}{cl}
2C_2(E)-\frac{(n-1)}{n}C_1^2(E)&= 
\frac{n_1n_2}{n}(\frac{\delta_1^2}{n_1^2}+\frac{\delta_2^2}{n_2^2}-
2\frac{\delta_1\delta_2}{n_1n_2})\\ 
&=-\alpha^2(\frac{n_1n_2}{n})\frac{V(d-1)!}{4\pi^2d}\ ,
\end{array}
\end{equation}
\noi where we have used $\alpha=\mu_1-\mu_2$\ in the last line.
\qed

\remarks 
\begin{enumerate}
\item The condition (\ref{eq:Fperp2}) can be applied to connections on complex 
bundles over symplectic manifolds, where $\omega$\ is then the 
symplectic form. It is thus tempting to view this as the 
definition a symplectic version of a stable Higgs extension, in 
much the same way that flat bundle provide the real versions of a 
stable Higgs bundles (under suitable restrictions on chern 
classes). However, as the above proof shows, the condition forces 
the Higgs extension to be a direct sum of polystable Higgs 
bundles, so no new phenomena emerge. It is also worth noting that, 
by (\ref{alpha}), the equation $F^{\perp}_H=-i\alpha\bT \omega $\ 
can apply only if $\alpha$\ is at the extreme lower bound of its 
range. 
\item In the case where $\Theta=0$, or equivalently $\nabla''=\dbare$,
Theorem \ref{thm:BogIneq} yields a Bogomolov Inequality for 
$\alpha$-stable extensions. This is equivalent to Theorem 3.11 in 
\cite{DUW}. Taking $\nabla''=\dbare$ in Theorem \ref{thm:BogEq} similarly 
yields a result for extensions of bundles. It provides the 
necessary and sufficient conditions under which equality can be 
attained in the Bogomolov inequality for an $\alpha$-stable 
extension. This result has not, as far as we are aware, previously 
appeared anywhere. 
 
\end{enumerate}
\section{Algebro-Geometric Description and \\
GIT Construction}\label{sect:GIT}
We now return to the algebraic setting and consider Higgs sheaves
and extensions of Higgs sheaves
as defined in Section \ref{sectn:Objects}. 
In \cite{DUW} Daskalopoulos, Uhlenbeck and Wentworth 
have constructed the moduli space of extensions of torsion free 
sheaves, following ideas of Simpson. In this Section we will show 
how, basically  the same proof of \cite{DUW}, also gives the 
moduli space of extensions of Higgs sheaves. The main modification 
required is to use sheaves of pure dimension, rather than torsion 
free sheaves. 

\par We will start by recalling Simpson's identification between Higgs 
sheaves on $X$ and sheaves on the cotangent bundle $T^*X$. Let $Z$ 
be the usual projective completion of the cotangent bundle 
$T^*X$, extending the projection $\pi:T^*X\lra X$ to a projective 
bundle $\overline\pi:Z\lra X$. Let $D=Z-T^*X$ be the divisor at 
infinity. Let $\SO_X(1)$ be an ample line bundle on 
$X$, and choose $b$ such that $\SO_Z(1):=\overline\pi^*\SO_X(b) 
\otimes \SO_Z(D)$ is an ample line bundle on $Z$. In \cite{S2} 
Simpson shows (cf. Lemma 6.8) that a Higgs sheaf 
$(\SE,\Theta)$ on $X$ is the same thing as a sheaf $\EE$ on $Z$ such
that $\supp(\EE)\cap D=\emptyset$. In fact, 
$\SE=\overline\pi_*\EE$, and the homomorphism $\Theta$ (with 
$\Theta \wedge \Theta =0$) is equivalent to giving the 
$\SO_{T^*X}$-module structure. This identification is also called 
the spectral cover construction. Denote $S=\supp(\EE)$, and 
consider the projection $\pi_S:S\lra X$. The fiber over a point 
$x\in X$ is a length $n=\rk(\SE)$, zero-dimensional subscheme of 
$T^*_x X=\Omega^1_x$, hence $\pi_S:S\lra X$ is an $n$-to-1 cover 
of $X$. If $X$ is a curve, then $S$ is the spectral curve studied 
in \cite{BNR}. The reason for this name is that if we restrict the 
Higgs field $\Theta$ to a point $x\in X$, we obtain an 
endomorphism of the fiber $E_x$ with values in 
$\Omega^1_x\cong \CC$ 
$$
\Theta_x: E_x \lra E_x \otimes \Omega^1_x,
$$
and hence the eigenvalues of $\Theta_x$ give a set of $n$ points 
(counted with multiplicity) of $T^*_x X$. This set is precisely 
the fiber of $S$ over $x\in X$. 

\par This identification between Higgs sheaves $(\SE,\Theta)$ on $X$ 
and torsion sheaves $\EE$ on $T^*_x X$ is compatible with 
morphisms, giving an equivalence of categories. The sheaf $\SE$ is 
torsion free if and only if $\EE$ is of pure dimension $d=\dim(X)$ 
(i.e., if $\EE$ is torsion free when restricted to its support and 
every irreducible component of its support has dimension $d$). 
Since 
$\SO_{T^*X}(1)=\pi^*\SO_X(b)$,
the Hilbert polynomials of $\EE$ and 
$\SE=\overline\pi_*\EE$ are related by 
$$
P(\EE,m)=P(\SE,bm)=:\wt P(\SE,m), 
$$
and hence $\EE$ is (semi)stable with respect to $\SO_X(1)$ if and 
only if $\SE$ is (semi)stable with respect to $\SO_Z(1)$ 
\cite[cor 6.9]{S2}. These correspondences between the Higgs sheaf 
and the sheaf of pure dimension are summarized in Table 2.

\begin{table}     
\begin{center}
\begin{tabular}{|c|c|c|}
\hline
\hline
& &  \\
& $(\cE,\Theta)$  Higgs sheaf on $X$ & $\EE$ sheaf on $T^*X$ \\
& & \\
\hline
& &  \\
\text{support}&\ $X$ &\ $S\subset T^*X$\ , spectral cover of $X$\\
& & \\
\hline
& &  \\
\text{Higgs structure} & 
$\Theta$ &
$\SO_{T^*X}$-module structure\\
& & \\
\hline
& &  \\
\text{sheaf type} & torsion free\ &\ of pure dimension $dim(X)$\\
& &  \\
\hline
& &  \\
\text{ample line bundle}\ &\ $\SO_X(1)$\ &\ 
$\SO_Z(1):=\overline\pi^*\SO_X(b) 
\otimes\SO_Z(D)$\\
& &  \\
\hline
& &  \\
\text{Hilbert polynomial}\ &\ $P(\SE,bm)$\ &\ 
$P(\EE,m)$\\
& & \\
\hline
& &  \\
\text{Gieseker stability}\ & w.r.t. $\SO_X(1)$\ & 
w.r.t. $\SO_Z(1)$\\ 
& &  \\
\hline
\hline
\end{tabular}
\begin{caption}
{}Algebro-Geometric Dictionary, giving the correspondence between 
Higgs sheaves on 
$X$ and sheaves of pure dimension on $T^*X\subset Z$ 
\end{caption}
\end{center}
\end{table}

Simpson then gives a method to construct the (projective) moduli 
space $\MP$ of semistable (with respect to $\SO_Z(1)$) sheaves 
with pure dimension on $Z$ and with Hilbert polynomial $\wt P$. 
Using the previous identification, plus the openness of the 
condition that $\supp(\EE)$ does not intersect 
$D$, one is thus able to identify $\MH$, the moduli space of 
semistable Higgs sheaves with Hilbert polynomial $P$, as an open 
subset of $\MP$. 

\par As in \cite{DUW}, instead of considering extensions, it is more 
convenient to take the equivalent point of view of considering 
\textit{quotient pairs} of Higgs sheaves. 
\begin{definition}\label{defn:quotPair}
A quotient pair of Higgs sheaves is a surjective morphism of Higgs 
sheaves 
$$
\CD\label{eq:quotPair}
(\SE,\Theta) @>{q}>> (\SF,\Psi) @>>> 0, 
\endCD$$
and it will be denoted by $q$ or by $(\SE, \Theta; \SF, \Psi)$. A 
morphism between quotient pairs of Higgs sheaves is a 
commutative diagram
\begin{equation}
\label{tomas1}
\CD
(\SE',\Theta') @>{q'}>>  (\SF',\Psi') @>>> 0 \\ 
      @V{f}VV              @V{g}VV    @. \\
(\SE,\Theta)   @>{q}>>   (\SF,\Psi)   @>>> 0 \\ 
\endCD
\end{equation}
\end{definition}

\noi {\bf Remark} Clearly, isomorphism classes of quotient pairs are 
the same thing as isomorphism classes of extensions. Indeed, using 
the notation of section \ref{sectn:Objects}, we take 
$(\SE_1,\Theta_1)=\ker q$, and 
$(\SE_2,\Theta_2)=(\SF,\Psi)$. We say that a quotient pair is stable
if the corresponding Higgs extension is stable.
A quotient pair $(\SE, \Theta; \SF, 
\Psi)$ is called torsion free if $\SE$ is a torsion free sheaf 
($\SF$ might have torsion). 
\begin{prop}[Jordan-H\"older filtration] If 
$(\SE, \Theta; \SF, \Psi)$ is a $\alpha$-Gieseker semistable torsion free
quotient pair, then there exists a filtration 
$$
\begin{array}{ccccccccccc}
(0,0) & = & (\SE_0,\Theta_0) & \subset & (\SE_1,\Theta_1) 
& \subset & \cdots & \subset & (\SE_l,\Theta_l) & = & (\SE,\Theta) \\
 &  & \Big\downarrow &  & \Big\downarrow 
&  &  &  & \Big\downarrow &  & \Big\downarrow \\
(0,0) & = & (\SF_0,\Psi_0) & \subset & (\SF_1,\Psi_1) 
& \subset & \cdots & \subset & (\SF_l,\Psi_l) & = & (\SF,\Psi) \\
 &  & \downarrow &  & \downarrow 
&  &  &  & \downarrow &  & \downarrow \\
 &  & 0 &  & 0 
&  &  &  & 0 &  & 0 \\
\end{array}
$$
such that $\SE_{i-1}$ is saturated in $\SE_i$ and the induced 
quotients 
$$
\overline{q}_i:(\SE_i/\SE_{i-1}, \overline{\Theta}_i) \lra
(\SF_i/\SF_{i-1}, \overline{\Psi}_i) 
$$
are $\alpha$-Gieseker stable and 
\begin{eqnarray*}
\frac{\deg(\SE_i/\SE_{i-1})-\alpha
\rk(\SF_i/\SF_{i-1})}{\rk(\SE_i/\SE_{i-1})} & = &
\frac{\deg(\SE)-\alpha \rk(\SF)}{\rk(\SE)}, \\
\frac{P(\SE_i/\SE_{i-1},m)}{\rk(\SE_i/\SE_{i-1})} & = &
\frac{P(\SE,m)}{\rk(\SE)} \quad \text{for all $m$, and} \\
\frac{P(\SF_i/\SF_{i-1},m)}{\rk(\SF_i/\SF_{i-1})} & = & 
\frac{P(\SF,m)}{\rk(\SF)} \quad \text{for all $m$.}
\end{eqnarray*}
Moreover, the direct sum of these quotient pairs, denoted 
$$
\gr(q)=\bigoplus_{i=1}^l  \overline{q}_i
$$
is unique up to isomorphism. 
\end{prop}
\pf
Analogous to \cite[prop 1.5.2]{HL} or \cite[prop 2.13]{DUW}. 
\qed

\noi {\bf Remark} Two quotient pairs $q$ and $q'$ are called 
S-equivalent if 
$\gr(q)\cong \gr(q')$. If $q$ is $\alpha$-Gieseker stable, then
$\gr(q)\cong q$.

\begin{thm} Fix Hilbert polynomials $P$ and $P''$.
There exists a quasi-projective scheme $\MHP$ whose points 
correspond to S-equivalence classes of quotient pairs of 
$\alpha$-Gieseker semistable torsion free Higgs sheaves with
the given Hilbert polynomials.
\end{thm}

\pf  The moduli space $\MTF$ of quotient pairs of torsion free sheaves 
has been constructed in \cite{DUW}, but since they use Simpson's 
method, their proof works not only for torsion free sheaves, but 
also for quotient pairs of sheaves of pure dimension. Let $\MPP$ 
be the moduli space of quotient pairs $\EE\lra \FF \to 0$ of 
sheaves on $Z$ with $\EE$ of pure dimension. Since the condition 
that $\supp(\EE)$ doesn't intersect $D$ is open, then using 
Simpson's identification we finally conclude that $\MHP$ is an 
open subset of $\MPP$. 

\par Now we will briefly recall the construction in \cite[section 
5]{DUW}, indicating what has to be changed to consider sheaves of 
pure dimension. For any coherent sheaf $\EE$ on $Z$, its Hilbert 
polynomial can be written as 
$$
\chi(\EE(m))=r\frac{m^d}{d!}+a \frac{m^{d-1}}{(d-1)!} +\cdots,
$$
where $d$ is the dimension of the support of $\EE$. Following 
Simpson \cite[p. 55]{S1}, we call $r$ the rank of $\EE$, and $a$ 
the degree of $\EE$ with respect to $\SO_Z(1)$. These definitions 
coincide with the usual definitions of rank and degree when $\EE$ 
is torsion free. If $\EE$ is a sheaf of pure dimension with 
support $S\subset Z$, then $r$ and $a$ are the rank and degree of 
$\EE$ when considered as a sheaf on its support $S$. 

\par Using these new definitions for rank and degree, the GIT 
construction in \cite{DUW} goes through for quotient pairs of pure 
dimension. First one proves that the set of semistable quotient 
pairs (with fixed Hilbert polynomials $\wt P$ and $\wt P''$) is 
bounded, and then there is an integer $K_0$ such that if $k\geq 
K_0$, for all semistable quotient pairs 
$q:\EE \lra \FF$ (with $\EE$ of pure dimension), $\EE(k)$ is generated
by global sections and $h^0(\EE(k))=\chi(\EE(k))=:N$. 

\par Let $V={\CC}^N$ be a fixed vector space of dimension $N$. Consider 
pairs $(q,\phi)$, where $q$ is a semistable quotient pair and 
$\phi:V\lra H^0(\EE(k))$ is an isomorphism. A pair $(q,\phi)$ is the
same thing as a commutative diagram 
\begin{equation}
\label{tomas2}
\begin{CD}
V\otimes \SO_Z @>{q_1}>> \EE(k) @>>> 0 \\ @|                            
@VV{q}V      \\ V\otimes \SO_Z @>{q_2}>> \FF(k) @>>> 0 \\ @.                            
@VVV         \\ @.                              0          \\ 
\end{CD}
\end{equation}
such that $q_1$ induces an isomorphism $V\cong H^0(\EE(k))$, hence 
for each pair $(q,\phi)$ we get a point $(q_1,q_2)$ in 
\begin{equation}
\label{tomas3}
\quot(V\otimes \SO_Z, \wt P_m) \times 
\quot(V\otimes \SO_Z, \wt P''_m)
\end{equation}
where $\quot(V\otimes \SO_Z, \wt P_m)$ (resp. $\quot(V\otimes 
\SO_Z, \wt P''_m)$) is Grothendieck's quotient 
scheme, parameterizing quotients of $V\otimes \SO_Z$ with Hilbert 
polynomial $\wt P_m(i)=\wt P(m+i)$ (resp. $\wt P''_m(i)=\wt 
P''(m+i)$). 

\par Let $\widehat Q_k$ be the closed subset of (3) where $\ker q_1 
\subset\ker q_2$ (i.e., $q_2$ factors through $q_1$), let 
$Q_k\subset \widehat Q_k$ be the subscheme where $\EE$ is of pure 
dimension, and let 
$\overline Q_k \subset \widehat Q_k$ be its closure. 
The projective scheme $\overline Q_k$ parameterizes commutative 
diagrams like (\ref{tomas2}). Now we have to get rid of the choice 
of isomorphism 
$\phi$. The group $\slv$ acts on (\ref{tomas3}) 
and hence on  $\overline Q_k$ (since this is invariant). From the 
point of view of pairs $(q,\phi)$, this action corresponds to 
$(q,\phi) \mapsto (q,g\circ \phi)$ for 
$g\in \slv$, so to get rid of the choice of the isomorphism $\phi$ we
only need to take the quotient by $\slv$. Note that it is enough 
to use $\slv$, and we don't need to use $\glv$, because scalar 
multiplication acts trivially on (\ref{tomas3}). 
This is done by taking the GIT quotient of $\overline Q_k$ by 
$\slv$, but to do this, first we have to linearize the action of $\slv$
on an ample line bundle on $\overline Q_k$. 
Following Grothendieck, by tensoring with $\SO_Z(j)$ for high 
enough $j$, and taking sections, we embed (\ref{tomas3}) (and 
hence $\overline Q_k$) into a product of Grassmanians 
$$
\Gr(V\otimes W, \wt P(k+j)) \times 
\Gr(V\otimes W, \wt P''(k+j)),
$$
where $W=H^0(\SO_Z(j))$. Using Pl\"ucker coordinates we get an 
embedding in 
\begin{equation}
\label{tomas4}
P=\PP\Big(\bigwedge {}^{\wt P(k+j)}(V\otimes W)^\vee\Big) \times 
\PP\Big(\bigwedge {}^{\wt P''(k+j)}(V\otimes W)^\vee\Big).
\end{equation}
The natural action of $\slv$ on (\ref{tomas4}) has a natural 
linearization on 
$\SO_P(r,s)$ for any $r$ and $s$, and by restriction we obtain a
linearization on the line bundle $\SO_P(r,s)|_{\overline Q_k}$ on 
$\overline Q_k$.  We choose $r$ and $s$ depending on
$\alpha$ as in \cite[p. 511]{DUW}.
Then one proves that GIT-semistable (resp. stable) points on 
$\overline Q_k$ correspond to $\alpha$-Gieseker semistable
(resp. stable) quotient pairs, and then the moduli space is 
obtained as the GIT quotient 
$$
\MPP= \overline Q_k /\!\!/ \slv.
$$
Finally one checks that points of $\MPP$ correspond to 
S-equivalence classes. 
\qed
\appendix
\numberwithin{equation}{subsection}
\section{Bott-Chern forms for Higgs Bundles}
\label{sectn:HBCclasses}
\subsection{ Introduction}\label{subs:intro}
In this Appendix adapt the computations of Bott and Chern (in 
their paper \cite{BC}) to construct Bott-Chern forms for Higgs 
Bundles. Keeping the notation of Section \ref{sectn:Eqtns}, 
\begin{itemize}
\item $\cE\lra X$\ is a rank $n$ holomorphic bundle with underlying smooth
complex bundle $E$ and holomorphic structure determined by an 
integrable partial connection $\dbare$ (as in \ref{eq:dbareDefn}),
\item A Higgs field on $E$\ is denoted by $\Theta$, and 
$\nabla''=\dbare+\Theta$\ is the Higgs operator. As in Definition 
\ref{defn:HiggsOpDesc}, a Higgs bundle on $X$\ is a pair 
$(E,\nabla'')$ in which $(\nabla'')^2=0$,
\end{itemize}
\begin{definition}\label{def:phi} Let $\phi$\ be any symmetric  $GL(n,\bC)$-invariant, 
$k$-linear function on $Mat_n$, the space of 
$n\times n$ matrices. We extend $\phi$\ to a k-linear map on
$Mat_n$-valued forms as follows: if $a_i\otimes\alpha_i\in
Mat_n\otimes\Omega^{p_i}(X)$, then 
\be\label{eq:phi}
\phi(a_1\otimes\alpha_1,\dots,a_k\otimes\alpha_k)=
\phi(a_1,\dots,a_n)\alpha_1\wedge\dots\wedge\alpha_k\ .
\end{equation}
\end{definition}
\noi Each $GL(n,\bC)$-invariant polynomial $\phi$ defines a 
characteristic class for $E$.  This class, denoted by $[\phi]\in 
H^{2k}(X,\bC)$, can be represented by the closed $2k$-form 
\be\label{eq:[phi]}
(\frac{i}{2\pi})^k\phi(F_D)\equiv(\frac{i}{2\pi})^k\phi(F_D,F_D,\dots,F_D)\ 
, 
\end{equation} 
\noindent where $D$\ is any 
$GL(n,\bC)$\ connection on $E$, and $F_D$\ is the 
$GL(n,\bC)$-valued 2-form which represents the 
curvature of $D$\ with respect to a local frame. Suppose now that 
$E$\ is the underlying smooth bundle of a holomorphic bundle $\cE=(E,\dbare)$.
Then any Hermitian bundle metric, say $H$, determines a unique 
Chern connection. Denoting the curvature of this connection by 
$F^D_H$, we thus get a representative 2k-form 
\be\label{eq:[phi(H)]}
(\frac{i}{2\pi})^k\phi(H)=(\frac{i}{2\pi})^k\phi(F^D_H)\ , 
\end{equation}
\noi corresponding to each metric. If $K$\ is any other metric then 
$\phi(K)$ and $\phi(H)$ must differ by a closed form since they 
represent the same class in cohomology. The Bott-Chern forms give 
a more refined measure of this difference between $\phi(K)$ and 
$\phi(H)$ for any pair of metrics. 

\par The essential ingredient in this construction is the Chern 
connection, which uses the defining structure of the holomorphic 
bundle (i.e. the operator $\dbare$) to associates a unique 
connection to each metric on $\cE=(E,\dbare)$. Suppose now that we 
add a Higgs field $\Theta$ to $\cE$ and, as outlined in \S 
\ref{sectn:Eqtns}, replace $\dbare$\ by the Higgs operator 
$\nabla''=\dbare+\Theta$. Each metric then produces a unique connection 
determined by the defining data of the Higgs bundle, i.e. 
determined by $\nabla''$ (or equivalently by $\dbare$ and 
$\Theta$). Given a $GL(n,\bC)$-invariant polynomial we can use 
these Higgs connections to associate to each metric, $H$, a Higgs 
representative for the corresponding characteristic class:
\begin{definition}\label{def:phiHiggs} Let $H$\ be a Hermitian metric on the
Higgs bundle $(E,\nabla'')$. Let $\nabla_H$\ be the corresponding 
Higgs connection, and let $\Fhiggs$\ be the curvature of this 
connection. Let $\phi$\ be any 
$GL(n,\bC)$-invariant, 
$k$-linear, symmetric function on $M_n$.  We define
\be\label{eq:phiHiggs}
\phi_{Higgs}(H)=\phi(\Fhiggs,\Fhiggs,\dots,\Fhiggs)\ .
\end{equation}
\end{definition}
The Higgs-Bott-Chern forms measure the difference between the 
closed forms $\phi_{Higgs}(H)$\ and 
$\phi_{Higgs}(K)$, for any two metrics $H$ and $K$. Our main 
result is as follows:
\begin{thm}\label{th:RHiggs}  Corresponding to each $GL(n,\bC)$-invariant, 
$k$-linear function $\phi$\ there is a function of 
pairs of metrics, $R_{Higgs}(H,K)$, such that: 
(i) $R_{Higgs}(H,K)$ takes its values in $\Omega ^{2k-1}(X,\bC)$, 
(ii)  $R_{Higgs}(H,K)$\ is well defined modulo $Im 
\overline{\partial} + Im \partial$, where $Im\overline{\partial}$ 
and $Im \partial$ denote the images 
$\overline{\partial}(\Omega ^{2k-1}(X,\bC))$ and 
$\partial(\Omega ^{2k-1}(X,\bC))$ in 
$\Omega ^{2k-1}(X,\bC)$, and 
(iii) 
\be\label{eq:phiHiggs(H,K)}
\phi_{Higgs}(H)-\phi_{Higgs}(K)=i\overline{\partial}\partial
R_{Higgs}(H,K)\ . 
\end{equation} 
\end{thm}
The forms $R_{Higgs}(H,K)$\ are the analogs for Higgs bundles of 
the Bott-Chern forms associated to pairs of metrics on a 
holomorphic bundle.  We will thus refer to these as 
\it Higgs Bott-Chern forms\rm .  Notice that unlike on holomorphic
bundles, for which the Bott-Chern forms take their values in 
$\Omega^{(k,k)}(X,\bC)$, the Higgs Bott-Chern forms need not
have holomorphic type $(k,k)$. This difference does not play any 
role in the proof of Theorem \ref{th:RHiggs}. Indeed the main 
ingredients in the proof are formally identical to those of 
Proposition 3.15 in \cite{BC}, the difference being that in place 
of the Chern connections used in \cite{BC}, here   we use Higgs 
connections. 
\subsection{Definition of $R_{Higgs}(H,K)$}\label{subs:defnR}

Fix $\phi$, a symmetric $GL(n,\bC)$-invariant k-linear function on 
$Mat_n$ as in Definition \ref{def:phi}. 

\par Notice that though $\phi$\ is symmetric, its extension to 
$Mat_n$-valued forms on $X$ is not in general symmetric because of the 
skew-symmetry of the wedge product on forms.  The symmetry will, 
however, be preserved if \it at most one\rm\ of the forms has odd 
degree. Since we will need them later, we record the following 
basic properties: 

\begin{lemma}\label{lem:dphi} Let $\phi$\ be any $GL(n,\bC)$-invariant,
$k$-linear function on $Mat_n$. For any matrix-valued forms 
$A_i=a_i\otimes\alpha_i\in Mat_n\otimes\Omega^{p_i}(X)$\ (for $i=1,\dots,k$), 
\be\label{eq:dphi}
d\phi(A_1,\dots,A_k)= \sum_j(-1)^{p_1+\dots+p_{j-1}} 
\phi(A_1,\dots,d(A_j),\dots,A_k)\ ,
\end{equation}
\noi If $B=b\otimes\beta\in Mat_n\Omega^q(X)$, then
\be\label{eq:symeqtn}
\sum_j(-1)^{p_{j+1}+\dots+p_{k}}\phi(A_1,\dots,[A_j,B],\dots,A_k)=0\ ,
\end{equation}
\noi where $[A_i,B]=[a_i,b]\alpha_i\wedge\beta$.
\end{lemma}
Given two metrics $H$\ and $K$\ we can pick a 1-parameter family 
of metrics, $H(t)$, such that $H(0)=H$ and $H(1)=K$, and so that 
it corresponds to a smooth path from $H$\ to $K$\ in the space of 
metrics.  We can compute derivatives with respect to the parameter
$t$ and thus define $L_t$\ by
\be\label{eq:Lt}
(L_ts,t)_{H(t)}=\frac{d}{dt}(s,t)_{H(t)}\ . 
\end{equation}
\begin{lemma}\label{Lemma:Lt} \cite{BC} Defined as above, $L_t$\ is a bundle
endomorphism, i.e. a global section in $\Omega^0(EndE)$. If $[H]$\ 
denotes the matrix representing $H$\ with respect to local frame 
$\{e_i\}$, then the matrix representing $L_t$\ is given by
\be\label{eq:LtMatrix}
 [L_t]=[H(t)]^{-1}[\dot{H(t)]}\ .
\end{equation}
\end{lemma}
\noindent Henceforth, where no confusion can arise, we drop the
square braces and denote the matrix representing $H$\ by $H$\ etc. 
Corresponding to the path of metrics $H(t)$ we get (cf. definition 
\ref{defn:HiggsConn}) a family 
\be\label{eq:nabla't}
\nabla'_t=D'_{H(t)} + \Theta^*_{H(t)}
\end{equation}
\noindent and thus a family of Higgs connections given by
\be\label{eq:nablat}
\nabla_t= \nabla''+\nabla'_t\ .
\end{equation}
\noindent Viewing the space of connections as an affine space, and 
identifying the tangent space at $\nabla_t$\ with 
$\Omega^1(X,EndE)$, we can compute the derivative with respect to $t$.
This yields an element $\dot{\nabla_t}\in\Omega^1(X,EndE)$. 
\begin{lemma}\label{lem:nabladot}\cite{BC}
\be\label{eq:nablatdot}
\frac{d}{dt}\nabla_t=\dot{\nabla_t}=\nabla'_t(L_t)\ ,
\end{equation}
\noindent where
\be
\nabla'_t(L_t)= \nabla'_t\circ L_t-L_t\circ\nabla'_t\ ,
\end{equation}
\noindent i.e. where $\nabla'_t(L_t)$\ is the contribution to the
covariant derivative $\nabla_t(L_t)$\ resulting from the 
decomposition of $\nabla_t$\ as $\nabla''+\nabla'_t$. 
\end{lemma}
We denote by $F_t$\ the curvature of the Higgs connection 
determined by $H(t)$, and define 
\be\label{eq:phi'}
\phi'_{Higgs}(F_t,L_t)=\sum_{j=1}^k\phi(F_t,\dots,
F_t,L_t,F_t\dots,F_t)\ , 
\end{equation}
\noi We compute 
\be
\begin{array}{cl}
\partial\phi'_{Higgs}(F_t,L_t)&=\sum_{i< j}
\sum_{j=1}^k\phi(F_t,\dots
\partial F_t,\dots,F_t,L_t,F_t\dots,F_t)+\nonumber\\ 
&\quad +\sum_{j=1}^k\phi(F_t,\dots,F_t,\partial 
L_t,F_t\dots,F_t)-\nonumber\\ 
&\quad -\sum_{i> j} \sum_{j=1}^k\phi(F_t,\dots,F_t,L_t,F_t\dots, 
\partial F_t,\dots,F_t)\ .
\end{array}
\end{equation}
\noi But by the Bianchi identities for Higgs connections, 
\be\label{eq:nabla'(Ft)}
\nabla'_t(F_t)=0=\partial F_t+[F_t,A_t]+[F_t,\Theta_t]\ ,
\end{equation}
\noindent where $\partial +A_t$\ is the $(1,0)$ part of the Chern
connection corresponding to $H(t)$.  Together with the invariance 
of $\phi$\ (cf. equation (\ref{eq:symeqtn}), and 
(\ref{eq:nablatdot}), this leads to the expression 
\be\label{eq: delphi'}
\begin{array}{cl}
\partial\phi'_{Higgs}(F_t,L_t)&=\sum_{j=1}^k\phi(F_t,\dots,F_t,\partial L_t
- [L_t,A_t]-[L_t,\Theta_t],F_t,\dots,F_t)\\ 
&=\sum_{j=1}^k\phi(F_t,\dots,F_t,\nabla'_t( L_t),F_t,\dots,F_t)\\ 
&=\phi'_{Higgs}(F_t,\dot{\nabla_t})\ .
\end{array}
\end{equation}
\noi But (cf. Proposition 2.18 in
\cite{BC}, or any standard discussion of the Chern-Weil homomorphism) 
$\int_0^1\phi'_{Higgs}(F_t,\dot{\nabla_t})dt$\ is precisely the
transgression term relating $\phi_{Higgs}(H)$\ and 
$\phi_{Higgs}(K)$, i.e.
\be\label{eq:trnsgrs}
\phi_{Higgs}(K)-\phi_{Higgs}(H)=d\left(\int_0^1\phi'_{Higgs}(F_t,\dot
{\nabla_t})dt\right)\ . 
\end{equation}
\noindent It thus follows from (\ref{eq: delphi'}) that
\be\label{eq:old25}\phi_{Higgs}(K)-\phi_{Higgs}(H)=
\overline{\partial}\partial\left(\int_0^1\phi'_{Higgs}(F_t,L_t)dt
\right)\ .
\end{equation}
\noindent We may therefore define
\begin{definition}\label{def:RHiggs}
 Given metrics $H$\ and $K$, and given a path 
$H(t)$ from $H$\ to $K$, set
\be\label{eq:RHiggs}
R_{Higgs}(H,K)=-i \int_0^1\phi'_{Higgs}(F_t,L_t)dt\ . 
\end{equation}
\end{definition}
\remark\ Notice in particular that (\ref{eq:old25}) implies that
$\overline{\partial}\partial R_{Higgs}(H,K)$\ is independent of
the path $H_t$\ joining $H$ and $K$. 
\subsection{Independence of the path $H(t)$}\label{subs:pathInd}
\noi In order to prove that $R_{Higgs}(H,K)$\ is well defined, i.e. 
is independent of the choice of path $H(t)$, we reformulate the 
definition in terms of a 1-form on $Met(E)$, the space of 
Hermitian metrics on $E$, and appeal to Stokes Theorem. Recall 
(cf. \cite{Ko}) that $Met(E)$ is a convex domain in an infinite 
dimensional vector space, and that the tangent space at any point 
$H\in Met(E)$ can be identified with hermitian sections of 
$End(E)$, i.e.
\be\label{eq:TMet}
T_HMet(E) = Herm_H(E)=\{u\in\Omega^0(EndE)\ |\ u^{*_H}=u\ \}\ .
\end{equation} 
\begin{definition}\label{def:theta} Let $U_H$\ be a tangent vector in
$T_HMet(E)$, and let $H(t)$\ be a path in $Met(E)$\ with $H(0)=H$\
and $\dot{H(0)}=U_H$. Define 
\be\label{eq:theta}
\theta_H(U_H)=\phi'_{Higgs}(\Fhiggs,L_0)\ ,
\end{equation}
\noindent where, as before, $L_t=H(t)^{-1}\dot{H(t)}$.
\end{definition}
\noi Given a curve $\gamma=H(t)$ which joins $H$ and $K$\  in
$Met(E)$, our definition of $R_{Higgs}(H,K)$\ thus becomes
\be\label{eq:Rtheta}
R_{Higgs}(H,K)=-i\int_{\gamma}\theta\ . 
\end{equation}
\noi  Expressed in this way, it becomes apparent that we can show the 
independence of the path $\gamma$\ by computing 
$d\theta$\ and applying Stokes Theorem. Suppose therefore that 
$U_H,V_H$\ are vectors in $T_HMet(H)$. Let $h(s,t)$\ be a smooth map from a
neighborhood of the origin in $\bR^2$\ to $Met(E)$, such that 
\be\label{eq:h(s,t)}
\begin{array}{cl}
 h(0,0) &= H \\ 
 h_*(\frac{\partial}{\partial s})&=
U_{h(s,t)}\ ,\\ 
 h_*(\frac{\partial }{\partial t})&= V_{h(s,t)}\ , 
\end{array}
\end{equation}
\noindent where $U_{h(s,t)}$\ and $V_{h(s,t)}$\ are vector fields
which extend $U_H$\ and $V_H$\ respectively.  Then 
\be\label{eq:dtheta}
\begin{array}{ll}
d\theta_H(U,V)&= h^*(d\theta)(\dds, \ddt)\\ 
&=\dds\theta(h_*\ddt) -\ddt\theta(h_*\dds)\\
&=U_H(\theta_H(V_H))-V_H(\theta_H(U_H))
\end{array}
\end{equation}
\begin{lemma}\label{lemma:TMet}
Under the identification of tangent spaces of 
$Met(E)$\ with hermitian sections of $End(E)$, as in (\ref{eq:TMet})
we get
\be\label{eq:dds}
\dds ( h^{-1}(s,t)V_{h(s,t)})|_{_{s=t=0}}= -U_HV_H
+H^{-1}\frac{\partial^2h}{\partial s\partial t}|_{_{s=t=0}} 
\end{equation}
\be\label{eq:ddsF}\dds F_{h(s,t)} =
\nabla''\nabla'_{h(s,t)}(h^{-1}(s,t)U_{h(s,t)})
\end{equation}
\end{lemma}
\noi We compute 
\be\label{eq:dtheta2}
\begin{array}{cl}
d\theta_H(U,V)= 
&\phi([H^{-1}V_H,H^{-1}U_H],\Fhiggs,\dots,\Fhiggs)-\\ &-
\sum_{j=2}^k\phi(H^{-1}U_H,\Fhiggs,\dots,\Fhiggs,\nabla''\nabla'_{H}
(H^{-1}V_{H}),\Fhiggs, 
\dots,\Fhiggs)+\\
&+\sum_{j=2}^k\phi(H^{-1}V_H,\Fhiggs,\dots,\Fhiggs,\nabla''\nabla'_{H}
(H^{-1}U_{H}),\Fhiggs, 
\dots,\Fhiggs)\ .
\end{array}
\end{equation}
\noi To simplify the notation, we set $u=H^{-1}U_H$\ and
$v=H^{-1}V_H$. The first term in (\ref{eq:dtheta2}) is then
\be\label{eq:phi[v,u]}
\begin{array}{cl}
\phi([v,u],\Fhiggs,\dots,\Fhiggs)&=-\sum_{j=2}^k\phi(v,\Fhiggs,\dots,\Fhiggs,[\Fhiggs,u],\Fhiggs,
\dots,\Fhiggs)\\ &=
-\sum_{j=2}^k\phi(v,\Fhiggs,\dots,\Fhiggs,\nabla''\nabla'_{H}(u),\Fhiggs, 
\dots,\Fhiggs)\\ &\quad
-\sum_{j=2}^k\phi(v,\Fhiggs,\dots,\Fhiggs,\nabla'_H\nabla''(u),\Fhiggs, 
\dots,\Fhiggs)\ ,
\end{array}
\end{equation}
\noi where the first equality follows by (\ref{eq:symeqtn}) and the second
equality follows from the fact that 
$$[\Fhiggs,u]=\Fhiggs(u)=\nabla''\nabla'_{H}(u)+\nabla'_H\nabla''(u)\ ,$$
\noi where the $\Fhiggs$\ in the expression $\Fhiggs(u)$\ refers to the
curvature of the induced connection on $End E$. Hence 
\be\label{eq:dtheta3}
\begin{array}{cl}
d\theta_H(U,V)&=-\sum_{j=2}^k\phi(u,\Fhiggs,\dots,\Fhiggs,\nabla''\nabla'_{H}(v), 
\Fhiggs, \dots,\Fhiggs)-\\ 
&\quad-\sum_{j=2}^k\phi(v,\Fhiggs,\dots,\Fhiggs,\nabla'_H\nabla''(u),\Fhiggs, 
\dots,\Fhiggs)\ .
\end{array}
\end{equation}
\begin{lemma}\label{lemma:dphi}
 For any connection $D$ on $E$, any
(symmetric), invariant k-linear function $\phi$, and any 
collection $A_i\in\Omega^{p_i}(End(E))$, for $i=1,\dots,k$, we 
have 
\be\label{eq:dphi2}
d\phi(A_1,\dots,A_k)= 
\sum_j(-1)^{p_1+\dots+p_{j-1}}\phi(A_1,\dots,DA_j,\dots,A_k)\ .
\end{equation}
\end{lemma}
\pf We fix a local frame for $E$ and write $D=d+A$, where
$A$\ is the connection 1-form. Thus $DA_j=dA_j+(-1)^{p_j}[A_j,A]$.
Using both parts of Lemma \ref{lem:dphi} we get 
\be\label{eq:dphi3}
\begin{array}{cl}
d\phi(A_1,\dots,A_k)&= 
\sum_j(-1)^{p_1+\dots+p_{j-1}}\phi(A_1,\dots,dA_j,\dots,A_k)\\ 
&=\sum_j(-1)^{p_1+\dots+p_{j-1}}\phi(A_1,\dots,DA_j,\dots,A_k)-\\
&\quad-\sum_j(-1)^{p_1+\dots+p_{j-1}+p_j}\phi(A_1,\dots,[A_j,A],\dots,A_k)\\ 
&=\sum_j(-1)^{p_1+\dots+p_{j-1}}\phi(A_1,\dots,DA_j,\dots,A_k) 
\end{array}
\end{equation}
\qed
\begin{cor}\label{cor:dbarphi}
 If $\nabla''=\dbare+\Theta$\ is the Higgs
operator, then 
\be\label{eq:dbarphi1}
\dbar\phi(A_1,\dots,A_k)=
\sum_j(-1)^{p_1+\dots+p_{j-1}}\phi(A_1,\dots,\nabla''A_j,\dots,A_k)\ ,
\end{equation}
\noi and if $\nabla'_H=D'_H+\Theta^*_H$, then
\be\label{eq:dbarphi2}
\partial\phi(A_1,\dots,A_k)=
\sum_j(-1)^{p_1+\dots+p_{j-1}}\phi(A_1,\dots,\nabla'_H
A_j,\dots,A_k)\ . 
\end{equation}
\end{cor}
\pf If we apply Lemma \ref{lemma:dphi} to the Chern connection
$\dbare+D'_H$, and decompose both side of (\ref{eq:dphi2}) according to
holomorphic type, we get 
\be\label{eq:dbarphi3}
 \dbar\phi(A_1,\dots,A_k)=
\sum_j(-1)^{p_1+\dots+p_{j-1}}\phi(A_1,\dots,\dbare
A_j,\dots,A_k) 
\end{equation}
\be\label{eq:dphi4}
 \partial\phi(A_1,\dots,A_k)=
\sum_j(-1)^{p_1+\dots+p_{j-1}}\phi(A_1,\dots,D'_H
A_j,\dots,A_k) 
\end{equation}
\noi But $\nabla''A_j=\dbare A_j +(-1)^{p_j}[A_j,\Theta]$. Equation
(\ref{eq:dbarphi3}) thus yields 
\be\label{eq:dbarphi4}
\begin{array}{cl}
\dbar\phi(A_1,\dots,A_k)&=\sum_j(-1)^{p_1+\dots+p_{j-1}}
\phi(A_1,\dots,\nabla''A_j,\dots,A_k)-\\
&-\sum_j(-1)^{p_1+\dots+p_{j-1}+p_j}\phi(A_1,\dots,[A_j,\Theta],\dots,A_k) 
\ .
\end{array}
\end{equation}
\noi The last summation in (\ref{eq:dbarphi4}) vanishes by (\ref{eq:symeqtn}) 
in Lemma \ref{lem:dphi}, i.e. by the invariance of $\phi$.  
Equation (\ref{eq:dbarphi2}) follows similarly from 
(\ref{eq:dphi4}), using the invariance of $\phi$ and
$\nabla'_H A_j= D'_H A_j+(-1)^{p_j}[A_j,\Theta^*_H]$. 
\qed

\bigskip

\noindent Using (\ref{eq:dbarphi1}) and (\ref{eq:dbarphi2}) of
Corollary \ref{cor:dbarphi}, the Bianchi identities 
(\ref{eq:BianId}), and Lemma 
\ref{lem:dphi} , we thus compute 
\be\label{eq:phi(nab''nab')}
\begin{array}{ll}
&\sum_{j=2}^k\phi(u,\Fhiggs,\dots,\Fhiggs,\nabla''
\nabla'_{H}(v),\Fhiggs,\dots,\Fhiggs)\\
&=-\sum_{j=2}^k\phi(\nabla''(u),\Fhiggs,\dots,\Fhiggs,\nabla'_{H}(v),\Fhiggs, 
\dots,\Fhiggs)-\dbar\alpha(u,v)
\end{array}
\end{equation}
\noi and
\be\label{eq:phi(nab'nab'')}
\begin{array}{ll}
&\sum_{j=2}^k\phi(v,\Fhiggs,\dots,\Fhiggs,\nabla'_H\nabla''(u),\Fhiggs,
\dots,\Fhiggs)\\ 
&=-\sum_{j=2}^k\phi(\nabla'_H(v),\Fhiggs,\dots,\Fhiggs,\nabla''(u),\Fhiggs, 
\dots,\Fhiggs)-\partial\beta(u,v)\ .
\end{array}
\end{equation}

\noi The forms $\alpha$\ and $\beta$\ are forms on $X$, given by 
\be\label{eq:alpha}
-\alpha(u,v)=\phi(u,\Fhiggs,\dots,\Fhiggs,\nabla'_H(v),\Fhiggs,\dots,\Fhiggs)\ 
\end{equation}
\noi and
\be\label{beta}-\beta(u,v)=
\phi(v,\Fhiggs,\dots,\Fhiggs,\nabla''(u),
\Fhiggs,\dots,\Fhiggs)\ . 
\end{equation}

\noi Furthermore, since $\nabla'_H(v)$\ and $\nabla''(u)$\ are
1-forms and $\Fhiggs$\ is a 2-form, it follows by the invariance 
of $\phi$\ (cf. the Remark after Definition \ref{def:phi}) that 
\be\label{eq:phi2}
\begin{array}{c}
\phi(\nabla'_H(v),\Fhiggs,\dots,\Fhiggs,\nabla''(u),\Fhiggs,
\dots,\Fhiggs)\\ +\phi(\nabla''(u),\Fhiggs,\dots,\Fhiggs,\nabla'_{H}(v),\Fhiggs,
\dots,\Fhiggs)=0\ .
\end{array}
\end{equation}

\noi Equation (\ref{eq:dtheta3}) thus reduces to
\be\label{eq: dtheta4}
d\theta_H(U,V)=\dbar\alpha(u,v)+\partial\beta(u,v) 
\end{equation}
\begin{lemma}\label{lemma:alpha+beta}
The expression 
$\dbar\alpha(u,v)+\partial\beta(u,v)$\ defines a 2-form on
$Met(E)$\ with values in $Im\dbar+Im\partial$
\end{lemma}
\pf Applying (\ref{eq:dbarphi2})in Corollary \ref{cor:dbarphi} to
$\phi(u,\Fhiggs,\dots,\Fhiggs,v,\Fhiggs,\dots,\Fhiggs)$\ gives
\be\label{eq:phi3}
\begin{array}{l}
\phi(u,\Fhiggs,\dots,\Fhiggs,\nabla'_H(v),\Fhiggs,\dots,\Fhiggs)\\ 
\qquad =-\phi(\nabla'_H(u),\Fhiggs,\dots,\Fhiggs,v,\Fhiggs,\dots,\Fhiggs)+ 
\partial\phi(u,\Fhiggs,\dots,\Fhiggs,v,\Fhiggs,\dots,\Fhiggs)\ ,
\end{array}
\end{equation}
\noi and hence
\be\label{eq:dbaralpha}
\dbar\alpha(u,v)=
\dbar\phi(\nabla'_H(u),\Fhiggs,\dots,\Fhiggs,v,\Fhiggs,\dots,\Fhiggs)
+\dbar\partial\phi(u,\Fhiggs,\dots,\Fhiggs,v,\Fhiggs,\dots,\Fhiggs) 
\end{equation}
\noi Similarly, applying (\ref{eq:dbarphi1}) to
$\phi(v,\Fhiggs,\dots,\Fhiggs,u,\Fhiggs,\dots,\Fhiggs)$\ gives
\be\label{eq:delbeta}
\begin{array}{c}
\partial\beta(u,v) =
\partial\phi(\nabla''(v),\Fhiggs,\dots,\Fhiggs,u,\Fhiggs,\dots,\Fhiggs)\\
+
\partial\dbar\phi(v,\Fhiggs,\dots,\Fhiggs,u,\Fhiggs,\dots,\Fhiggs)\ .
\end{array}
\end{equation}
\noi Notice that in each occurrence of $\phi$ in (\ref{eq:dbaralpha}) and 
(\ref{eq:delbeta}), the arguments include at most one form of odd 
degree. By the remark after Definition \ref{def:phi} the 
expressions are thus symmetric functions of their arguments. 
Recall also that 
$\dbar\partial+\partial\dbar=0$. Combining (\ref{eq:dbaralpha}) and 
(\ref{eq:delbeta}) thus yields 
\be\label{eq:alpha+beta2}
\begin{array}{cl}
\dbar\alpha(u,v)+\partial\beta(u,v)&= \dbar\phi(v,
,\Fhiggs,\dots,\Fhiggs,\nabla'_H(u),\Fhiggs,\dots,\Fhiggs)\\ 
&\quad+\partial\phi(u,\Fhiggs,\dots,\Fhiggs,\nabla''(v),\Fhiggs,\dots,\Fhiggs)\\
&=-(\dbar\alpha(v,u)+\partial\beta(v,u))\ .
\end{array}
\end{equation}
\qed

\noi We can now prove
\begin{prop}\label{prop:pathInd}
Up to terms in $Im\partial+Im\dbar$, 
$R_{Higgs}(H,K)$\ is independent of the path $H(t)$\ used to
compute it in Definition \ref{eq:RHiggs}. Thus the map 
\be 
H\longmapsto R_{Higgs}(H,K) 
\end{equation}
\noindent gives a well defined map from $Met(E)$\ (the space of
metrics) to $\Omega^{k}(X,\bC )/Im\partial+Im\overline{\partial}$. 
\end{prop}
\pf Let $\gamma_1\ ,\ \gamma_2$\ be any
two paths  from $H$\ to $K$ in $Met$. Then $\gamma_1-\gamma_2$\ 
bounds a disk, say $\Gamma$, and Stokes Theorem implies 
\be
\int_{\gamma_1}\theta-\int_{\gamma_2}\theta=\int_{\Gamma}d\theta
=\int_{\Gamma}(\dbar\alpha + \partial\beta) \ . 
\end{equation}
\qed

\bigskip

\par The rest of Theorem \ref{th:RHiggs} now follows from the 
definition of $R_{Higgs}$. 

\remark It follows from the definition of $R_{Higgs}$\ that
if $H(t)$\ is a smooth 1-parameter family of metrics, then 
\be\label{RHiggsdot}
\frac{d}{dt}R_{Higgs}(H(t),K)=-ik\phi(L_t,F_t,\dots,F_t)\ ,
\end{equation}
\noi where $L_t$\ is as in (\ref{eq:Lt}) and $F_t$\ is the curvature of 
the Higgs connection corresponding to $H(t)$. 
 
\subsection{Two Special Cases}\label{subs:specialcases}
\bigskip
\subsubsection{Case 1}  
If $k=1$\ and $\phi(A)=Tr(A)$, then 
\be
\phi'(F_t,L_t)=\phi(L_t)=Tr(\dot{H(t)}H(t)^{-1})\ .
\end{equation}
\noi Thus, denoting the corresponding function $R_{Higgs}$\ by 
$R^{(1)}_{Higgs}$, we get
\be
R^{(1)}_{Higgs}(H,K)=-i\int_0^1Tr(\dot{H(t)}H(t)^{-1})dt\ . 
\end{equation}
\noindent Notice that this is the same as the corresponding
Bott-Chern form defined on a holomorphic bundle. In both cases 
(i.e. with or without the extra Higgs bundle structure) we get 
\be
R^{(1)}_{Higgs}(H,K)=-i\ln HK^{-1}\ , 
\end{equation}
\noi  which is manifestly independent of the path from $H$ to $K$.
\bigskip
\subsubsection{Case 2} 
If $k=2$ and 
$\phi(A_1,A_2)=-\frac{1}{2}Tr(A_1A_2+A_2A_1)$, then
\be
\phi'(F_t,L_t)=\phi(F_t,L_t)=-Tr(F_tL_t)\ ,
\end{equation}
$$R^{(2)}_{Higgs}(H,K)=i\int_0^1Tr(F_tL_t)dt\ .$$
\noindent The functional defined by Simpson in \cite{S3} is
\be\label{eq:MSimpson}
M_{S}(H,K)=\int_XR^{(2)}_{Higgs}(H,K)\wedge\omega^{d-1}\ . 
\end{equation}
\noindent This is the Higgs analog
of the function defined by Donaldson in \cite{Do1}, which is given 
by the same formula but with the Bott-Chern form $R^{(2)}(H,K)$ in 
place of the Higgs Bott-Chern form $R^{(2)}_{Higgs}(H,K)$. 
\begin{prop}\label{prop: MHiggsderivs}
Take $H(t)=Ke^ts$, with $s=s^{*_K}$. Then 
\be
\frac{d}{dt}M_{S}(H(t),K)=-2i\int_X\phi'(F_t,s)\wedge\omega^{d-1}
=2i\int_X Tr(F_ts)\wedge\omega^{d-1} 
\end{equation}
\be 
\frac{d^2}{dt^2}M_{S}(H(t),K)|_{t=0}=  |\nabla''(s)|_K^2
\end{equation}
\end{prop}
\pf The formulae for $\frac{d}{dt}M_{S}$\ follow 
directly from (\ref{RHiggsdot}). Using this result, plus the fact 
that (cf. (\ref{eq:ddsF})) 
$\frac{dF_t}{dt}=\nabla''\nabla'_t(s)$, we get 
\be
\begin{array}{cl}
\frac{d^2}{dt^2}M_{S}(H(t),K)|_{t=0}&=
2i\int_X Tr(\nabla''\nabla'_K(s)s)\wedge\omega^{d-1}\\ 
&=-2i \int_X Tr(\nabla''(s)\wedge \nabla'_K(s))\wedge\omega^{d-1}\\ 
&=2 \int_X |\nabla''(s)|^2_K\wedge\omega^{d-1}\ .
\end{array}
\end{equation}
\noi The second equality follows by (\ref{eq:phi(nab''nab')}). The third follows by 
Lemma 3.1(b) in \cite{S3}. 
\qed


\begin{thebibliography}{DUW}
\bibitem[BC]{BC}{R. Bott and S.S. Chern, }
\textit{Hermitian vector bundles and the equidistribution of the zeroes of 
their holomorphic sections, } Acta. Math. \textbf{114}, 71--112 
(1965). 
\bibitem[BNR]{BNR}{A. Beauville, M.S. Narasimhan and S. Ramanan, }
\textit{Spectral curves and the generalised theta divisor, }
J. reine angew. Math. \textbf{398}, 169--179 (1989). 
\bibitem[B1]{B1}{S. Bradlow, }
\textit{Special metrics and stability for holomorphic bundles with 
global sections, }J.D.G.\textbf{33}, 169--213 (1991) 
\bibitem[BGP]{BGPcag}{S. Bradlow and O. Garcia-Prada, }
\textit{Higher cohomology triples and holomorphic extensions, } Comm. 
Anal. Geom. \textbf{3}, 421--463 (1995). 
 
\bibitem[DUW]{DUW}{G. Daskalopoulos, K. Uhlenbeck, and 
R. Wentworth, } 
\textit{Moduli of extensions of holomorphic bundles on K\"ahler 
manifolds, } Comm. Anal. Geom. \textbf{3}, 479--522 (1995). 
\bibitem[Do1]{Do1}{S. Donaldson, } 
\textit{Anti-self-dual Yang-Mills connections on a complex 
algebraic surface and stable vector bundles, }Proc. Lond. Math. 
Soc.  \textbf{50}, 1--26 (1985). 
\bibitem[Do2]{Do2}{S. Donaldson, } 
\textit{Infinite determinants, stable bundles and curvature, }Duke Math. J.  
\textbf{54}, 231--247 (1987).
 
\bibitem[H]{H}{N. Hitchin, }
\textit{The self-duality equations on a Riemann 
surface} Proc. Lond. Math. Soc.\textbf{55}, 59-126, (1987)
\bibitem[HL]{HL}{D. Huybrechts and M. Lehn, }
\textit{The geometry of moduli spaces of sheaves, }
Aspects of Mathematics E31, Vieweg, Braunschweig/Wiesbaden 1997. 
\bibitem[Ko]{Ko}{S. Kobayashi, }
\textit{Differential Geometry of Complex Vector Bundles, 
}Princeton, 1987.
\bibitem[S3]{S3}{C. Simpson, }
\textit{Constructing variations of Hidge structure using Yang-Mills theory
and applications to uniformization, } JAMS 
\textbf{1}, 867-918 (1988).
\bibitem[S4]{S4}{C. Simpson, }
\textit{Higgs bundles and local systems, } Publ. Math. I.H.E.S. \textbf{75}, 
5--95 (1992). 
\bibitem[S1]{S1}{C. Simpson, }
\textit{Moduli of representations of the fundamental group of a smooth
projective variety I, } Publ. Math. I.H.E.S. \textbf{79}, 47--129 
(1994). 
\bibitem[S2]{S2}{C. Simpson, }
\textit{Moduli of representations of the fundamental group of a smooth
projective variety II, } Publ. Math. I.H.E.S. \textbf{80}, 5--79 
(1995). 
\bibitem[UY]{UY}{K. Uhlenbeck and S.T. Yau, }
\textit{On the existence of Hermitian-Yang-Mills connections in 
stable vector bundles} Comm. Pure and Appl. Math.\textbf{39-S}, 
257--293 (1986). 
\end{thebibliography}
\end{document}